\documentclass[letterpage, 11pt, notitlepage]{article}
\usepackage[margin=1.0in]{geometry}

\usepackage[round]{natbib}

\usepackage[utf8]{inputenc} 
\usepackage[T1]{fontenc}    
\usepackage{hyperref}       
\usepackage{url}            
\usepackage{booktabs}       
\usepackage{amsfonts}       
\usepackage{nicefrac}       
\usepackage{microtype}      
\usepackage{bbold}
\usepackage{capt-of}

\usepackage[center]{caption}

\usepackage{amsmath, amsthm, amssymb, graphicx, cite,float, epsfig,epstopdf,color,soul,mathabx}
\usepackage{tabularx}
\usepackage[ruled,vlined]{algorithm2e}
\allowdisplaybreaks
\usepackage{color}
\usepackage{enumerate}
\usepackage[shortlabels]{enumitem}
\usepackage{multicol}
\usepackage{empheq}
\usepackage{caption} 

\newcommand{\Rset}{\mathbb{R}}

\newcommand{\Acal}{{\cal A}}

\newcommand{\Fcal}{{\cal F}}

\newcommand{\Ocal}{{\cal O}}
\newcommand{\Pcal}{{\cal P}}

\newcommand{\Scal}{{\cal S}}

\newcommand{\Xcal}{{\cal X}}

\newcommand{\1}{{\mathbf{1}}}

\newcommand{\argmin}{\mathop{\rm argmin}}
\newcommand{\argmax}{\mathop{\rm argmax}}

\newtheorem{prop}{Proposition}
\newtheorem{lem}{Lemma}
\newtheorem{thm}{Theorem}

\newtheorem{definition}{Definition}

\newtheorem{assump}{Assumption}

\usepackage{tikz}

\usepackage[toc,page,header]{appendix}
\usepackage{minitoc}

\title{Finite-Time Complexity of Online Primal-Dual Natural Actor-Critic Algorithm for Constrained Markov Decision Processes}
\author{Sihan Zeng\thanks{School of Electrical and Computer Engineering, Georgia Institute of Technology, Atlanta, GA
.}
\and Thinh T. Doan\thanks{Bradley Department of Electrical and Computer Engineering, Virginia Tech, Blacksburg, VA.}
\and Justin Romberg\footnotemark[1]}

\begin{document}
\maketitle


\begin{abstract}
We consider a discounted cost constrained Markov decision process (CMDP) policy optimization problem, in which an agent seeks to maximize a discounted cumulative reward subject to a number of constraints on discounted cumulative utilities. 
To solve this constrained optimization program, we study an online actor-critic variant of a classic primal-dual method where the gradients of both the primal and dual functions are estimated using samples from a single trajectory generated by the underlying time-varying Markov processes. This online primal-dual natural actor-critic algorithm maintains and iteratively updates three variables: a dual variable (or Lagrangian multiplier), a primal variable (or actor), and a critic variable used to estimate the gradients of both primal and dual variables. These variables are updated simultaneously but on different time scales (using different step sizes) and they are all intertwined with each other. Our main contribution is to derive a finite-time analysis for the convergence of this algorithm to the global optimum of a CMDP problem. Specifically, we show that with a proper choice of step sizes the optimality gap and constraint violation converge to zero in expectation at a rate $\Ocal(1/K^{1/6})$, where $K$ is the number of iterations. To our knowledge, this paper is the first to study the finite-time complexity of an online primal-dual actor-critic method for solving a CMDP problem. We also validate the effectiveness of this algorithm through numerical simulations.
\end{abstract}

\section{Introduction}

The Markov decision process (MDP) is a mathematical framework for sequential decision making problems, where an agent aims to maximize the cumulative reward collected from an environment by taking a sequence of actions \citep{sutton2018reinforcement,Puterman_book_1994}. In many safety-centric applications in robotics \citep{koppejan2011neuroevolutionary,ono2015chance,hakobyan2019risk}, power systems \citep{vu2020safe}, and asset management \citep{abe2010optimizing}, besides maximizing the cumulative reward, it is equally, if not more, important for the agent to respect safety constraints and/or adopt reasonable behavior in the environment.

One way that safe behavior can be enforced in sequential decision making problems is to introduce constraints in the MDP.  In a  constrained Markov decision process (CMDP) problem, an agent maximizes its expected cumulative reward subject to constraints on its expected cumulative utilities. If the transition probability kernel of the environment is known, CMDP problems can be solved by dynamic programming \citep{altman1993asymptotic,altman1999constrained,piunovskiy2006dynamic}. When the transition probabilities are unknown but the agent can collect samples by interacting with the environment, 
reinforcement learning (RL) algorithms such as Q learning and actor-critic methods have been used to solve the CMDP problems.  Existing theoretical results in this domain have established the asymptotic convergence of these algorithms \citep{borkar2005actor,lakshmanan2012novel,bhatnagar2012online}.


In this paper, motivated by recent advances in the theoretical understanding of policy gradient methods \citep{ding2020natural,agarwal2020optimality,mei2020global}, we study the finite-time complexity of the so-called online primal-dual natural actor-critic algorithm for solving CMDP problems. This primal-dual method performs gradient descent on the dual variable in the Lagrangian of the underlying optimization program while a natural actor-critic algorithm ascends on the primal variable (the policy).  
We are interested in the online setting, where the policy maintained by the actor and the corresponding value function maintained by the critic are updated simultaneously using data sampled from a single trajectory generated by the time-varying policy at the actor.  Our goal is to derive the finite-time complexity of this online primal-dual natural actor-critic algorithm.

Compared with the ``batch'' (or off-line) setting where the critic update requires sampling from multiple re-started trajectories \citep{qiu2021finite}, our observation model of a single trajectory is more realistic and sample-efficient.  Its analysis, however, introduces several new technical challenges.  The batch actor-critic algorithm decouples the critic from the actor by allowing multiple updates of the critic using independently generated trajectories at every iteration of the actor update.  In the online algorithm, the actor and critic update simultaneously with every single sample.  These updates are coupled with one another, but occur on two different ``time scales'', meaning that dramatically different step sizes are used for their respective updates.  Together with the dual variable, there are three intertwined updates in each iteration, resulting in a multiple-time-scale method. Analyzing the performance of this method, especially its finite-time complexity, is complicated; we detail some the technical challenges in Section~\ref{sec:mainresults:challenges}.


\subsection{Main Contributions}

The main contribution of this paper is to provide a finite-time analysis of the primal-dual natural actor-critic algorithm for solving the CMDP problem.  The algorithm consists of multiple updates at different time scales performed online as data is continuously sampled from a single trajectory. Specifically, we show that the objective values of the policy iterates converge to the globally optimal value while the magnitude of their constraint violations converge to zero, both at a rate of $\Ocal(1/K^{1/6})$, where $K$ is the number of iterations. 
This translates to a sample complexity of $\Ocal(1/\delta^6)$ for the error to reach a precision of $\delta$. To the best of our knowledge, this is the first known finite-time and finite-sample analysis of the online primal-dual actor-critic algorithm for solving CMDP problems. Finally, we illustrate and compare the behavior of this algorithm with other methods through numerical simulations.

\subsection{Related Works}

There has been a surging interest in the machine learning community in developing safe and robust RL algorithms, where CMDP is a common mathematical framework for modelling safety constraints \citep{feyzabadi2014risk,brunke2021safe}. Different variants of primal-dual methods have been proposed to solve the CMDP problem and have exhibited strong empirical performance in many real-world applications \citep{achiam2017constrained,liang2018accelerated,ray2019benchmarking, spooner2020natural}. 

From the theoretical perspective, local and asymptotic convergence of various primal-dual algorithms have been studied in \citet{chow2017risk,chow2018lyapunov,tessler2018reward, yu2019convergent}. Despite the objective being non-concave and the constraint set non-convex, it has been shown that the duality gap of a CMDP problem is zero \citep{altman1999constrained,paternain2019constrained}. Based on this property, the authors in \citet{ding2020natural} analyze a primal-dual natural policy gradient algorithm for a CMDP with a single constraint and present a finite-time convergence to global optimality in the tabular policy setting and with a general smooth policy class. Another related work \citep{ding2021provably} studies an optimism-driven primal-dual algorithm for CMDP with a linear transition kernel and also provides a global convergence guarantee.

These works make the critical assumption of an oracle that returns an exact policy gradient at every iteration (or alternatively that the gradient is computed off-line using multiple samples generated by the current policy at every iteration). This setting is obviously not applicable to problems where data is generated in online fashion (e.g., in robotics and self-driving cars), and 
online algorithms are needed that update the variable iterates using gradients estimated by only one single data point at any iteration. The online setting, however, is more challenging to analyze, as we will see in Section~\ref{sec:mainresults:challenges}. 


There are also prior works in the literature on the convergence of online primal-dual actor-critic algorithms for CMDP policy optimization. For example, \citet{borkar2005actor} studies the asymptotic convergence of this method for the tabular setting, while \citet{bhatnagar2012online} considers convergence with linear function approximation. To our best knowledge, there does not exist any work that derives the finite-time complexity (or convergence rates) of the online primal-dual actor critic algorithm. Our focus, is, therefore, to fill in this gap.

Our approach is inspired by recent works on the finite-time analysis of policy gradient and actor-critic algorithms. In particular, the authors in \citet{agarwal2020optimality} show that the unconstrained MDP problem, despite being non-convex, has a ``gradient domination'' property, which essentially means that every stationary point is globally optimal. Using this property, different variants of policy gradient methods have been shown to return a global optimal policy \citep{agarwal2020optimality,mei2020global}. 
For the linear quadratic regulator (LQR) problem, the stronger Polyak-Lojasiewicz condition has been shown to hold and has been
used to provide convergence guarantees for the policy gradient method 
\citep{fazel2018global,yang2019global,zeng2021two}. 
Online actor-critic methods for unconstrained MDPs have also been studied in both the tabular \citep{khodadadian2021finite} and linear function approximation \citep{wu2020finite} settings.

Finally, we note that safe RL has also been studied in different frameworks. Examples of non-CMDP methods include works that use transfer learning or prior knowledge \citep{taylor2009transfer,anwar2020autonomous,garcia2020teaching}, and/or consider the minimax formulation under environment uncertainty \citep{pinto2017robust,ren2020improving}. For more discussion along this direction, we refer interested readers to the survey papers on safe and robust RL \citep{garcia2015comprehensive,kim2020safe}.
We also note that the CMDP policy optimization problem can be treated by algorithms beyond the primal-dual type. For example, \citet{xu2021crpo} proposes a primal-only approach to this problem with a provable convergence guarantee.

\section{Constrained Markov Decision Processes}
We consider a constrained Markov decision process with $N$ constraints specified by the 5-tuple $(\Scal,\Acal,\Pcal,\gamma, \{r_i\}_{i=0}^N)$, where $\Scal$ is the state space, $\Acal$ is the action space, $\Pcal$ is the set of transition probabilities represented the dynamics of the states, $\gamma\in(0,1)$ is the discount factor, $r_0:\Scal\times\Acal\rightarrow[0,1]$ is the reward function, and  $r_1,...,r_N:\Scal\times\Acal\rightarrow[0,1]$ are the utility functions that are used to define the constraints. We assume that the state and action space are finite.
%

Given a policy $\pi$, we define its Q function and value function associated with the reward/utility as the expectation of the discounted cumulative reward/utility
\begin{align*}
    &Q_{i}^{\pi}(s,a)\triangleq\mathbb{E}_{\pi}\Big[\sum_{t=0}^{\infty} \gamma^k r_i\left(s_k, a_k\right) \mid s_{0}=s, a_{0}=a\Big],\notag\\
    &V_{i}^{\pi}(s)\triangleq\mathbb{E}_{\pi}\Big[\sum_{k=0}^{\infty} \gamma^k r_i\left(s_k, a_k\right) \mid s_{0}=s\Big],\,\forall i=0,\cdots,N.
\end{align*}
We also define the advantage function as
\begin{align}
    A_{i}^{\pi}(s,a)\triangleq Q_{i}^{\pi}(s,a)-V_{i}^{\pi}(s),\,\,\forall i=0,\cdots,N.
\end{align}
With an abuse of notation, we use $V_{i}^{\pi}(\rho)$ to denote the expectation of the value function under policy $\pi$, reward/utility function $r_i$ with respect to an initial distribution $\rho$
\begin{align}
    V_{i}^{\pi}(\rho)\triangleq\mathbb{E}_{s_{0} \sim \rho}\left[V_{i}^{\pi}\left(s_{0}\right)\right],\,\,\forall i=0,\cdots,N.
\end{align}
The objective of the CMDP problem is to find a policy that maximizes the expected discounted cumulative reward while taking into account the constraints on the expected discounted cumulative utilities, i.e., we aim to solve 
\begin{align}
    \pi^{\star}=&\argmax_{\pi} \quad V_{0}^{\pi}(\rho)\notag\\
    &\text{subj. to} \quad V_{i}^{\pi}(\rho) \geq b_i\quad\forall i = 1,2,...,N,
    \label{eq:objective}
\end{align}
for prescribed constants $b_1, b_2,..., b_N$. 
We denote $V_g^{\pi}(\rho) = [V_1^{\pi}(\rho),..., V_N^{\pi}(\rho)]^{\top}$ and $b=[b_1,...,b_N]^{\top}$. In addition, we consider the following classic Slater's condition for constrained optimization problems. 

\begin{assump}[Slater's Condition]\label{assump:Slater}
There exists a constant $0<\xi\leq 1$ and a policy $\pi$ such that $V_i^{\pi}(\rho)\geq\xi$ for all $i=1,\cdots,N$.
\end{assump}
We denote the Lagrangian of problem \eqref{eq:objective} as
\begin{align}
 V_{L}^{\pi,\lambda}(\rho)\triangleq V_{0}^{\pi}(\rho)+\lambda^{\top}\left(V_{g}^{\pi}(\rho)-b\right),
    \label{eq:Lagrangian}
\end{align}
where $\lambda$ is the dual variable (or Lagrangian multiplier). The dual function $V_{D}^{\lambda}$ 
is defined as
\begin{align}
    V_D^{\lambda}(\rho)=\max_{\pi} V_{L}^{\pi,\lambda}(\rho),
    \label{eq:def_dualfunc}
\end{align}
and the dual problem is
\begin{align}
    \lambda^{\star}=
    \argmin_{\lambda\in\mathbb{R}_{+}^{N}} V_D^{\lambda}(\rho).
    \label{eq:def_dualprob}
\end{align}
Although neither the objective nor constrain set of CMDP is convex, the strong duality holds for \eqref{eq:objective} under Assumption \ref{assump:Slater} \citep[Theorem 3.6]{altman1999constrained}
\begin{align}
    V_{D}^{\lambda^{\star}}(\rho)=V_{0}^{\pi^{\star}}(\rho),\label{lem:strong_duality}
\end{align}
where $\pi^{\star}$ and $\lambda^{\star}$ are the (not necessarily unique) optimal solutions to \eqref{eq:objective} and \eqref{eq:def_dualprob}, respectively. This also implies that solving \eqref{eq:objective} is equivalent to finding a saddle point $(\pi^{\star},\lambda^{\star})$ of \eqref{eq:Lagrangian}.
Finally, under Slater's condition we show the boundedness of $\lambda^{\star}$ in the following lemma, which is an extension of \citet[Lemma 1]{ding2020natural} to the case of multiple constraints. We present its proof in the appendix. 
\begin{lem}\label{lem:bounded_lambdastar}
Under Assumption \ref{assump:Slater}, we have
\begin{align*}
    \|\lambda^\star\|_{\infty}\leq \frac{2}{\xi(1-\gamma)}\cdot
\end{align*}
\end{lem}

\section{Online Primal-Dual Natural Actor-Critic Algorithm}
In this section, we present an online actor-critic variant of the classic primal-dual methods for finding a saddle point of \eqref{eq:Lagrangian}, which is formally stated in Algorithm~\ref{Alg:AC}. We consider the tabular setting for both the value functions and the policy. Specifically, we maintain a table $\widehat{Q}_i\in\mathbb{R}^{|\Scal|\times|\Acal|}$ (critic) as an estimate of the Q function under the reward or utility function $r_i$ and a policy parameter table $\theta\in\mathbb{R}^{|\Scal|\times|\Acal|}$ (actor) which parametrizes the policy through the softmax function
\begin{align*}
    \pi_{\theta}(a \mid s)=\frac{\exp \left(\theta_{s, a}\right)}{\sum_{a^{\prime} \in \Acal} \exp \left(\theta_{s, a^{\prime}}\right)},\quad \text {for all} \quad \theta \in \mathbb{R}^{|\Scal||\Acal|}.
\end{align*}
Implementing the primal-dual method requires estimating the gradients of the Lagrangian $V_L^{\pi_{\theta},\lambda}$ with respect to the primal and dual variables, $\theta$  and  $\lambda$, respectively. It has been shown in \citet{ding2020natural} that $\widetilde{\nabla}_{\theta}V_L^{\pi_{\theta},\lambda}$, the natural gradient of $V_L^{\pi_{\theta},\lambda}$ with respect to $\theta$ and $\nabla_{\lambda}V_L^{\pi_{\theta},\lambda}$, the gradient of $V_L^{\pi_{\theta},\lambda}$ with respect to $\lambda$ are
\begin{align*}
    \widetilde{\nabla}_{\theta_{s,a}}V_L^{\pi_{\theta},\lambda}(\rho)&=\frac{1}{1-\gamma}\left(A_0^{\pi_{\theta}}(s,a)+\sum_{i=1}^{N}\lambda_i A_i^{\pi_{\theta}}(s,a)\right),\\
    \nabla_{\lambda_i}V_L^{\pi_{\theta},\lambda}(\rho)&=V_i^{\pi_{\theta}}(\rho)-b_i\notag\\
    &=\sum_{s,a}\rho(s) \pi_{\theta}(a\mid s) Q_i^{\pi_{\theta}}(s,a)-b_i.
\end{align*}
Evaluating the (natural) gradients requires knowing exactly the value functions of the policy. 
In our online algorithm, we cannot compute the value functions exactly and instead use their approximates estimated by the critic. 
This results in the dual variable update in \eqref{Alg:AC:dual_update} and the following update for $\theta$
\begin{align}
    \theta_{k+1}(s_k,a_k)=\theta_{k}(s_k,a_k)+\alpha \widehat{Q}_{L,k}(s, a),\label{Alg:AC:actor_update}
\end{align}
where we define
\begin{align*}
    &\widehat{Q}_{L,k}(s, a)\triangleq \widehat{Q}_{0,k}(s, a)+\lambda_k^{\top} \widehat{Q}_{g,k}(s, a)\in\mathbb{R},\\
    &\widehat{Q}_{g,k}(s,a) \triangleq [\widehat{Q}_{1,k}(s,a),...,\widehat{Q}_{N,k}(s,a)]^{\top}\in\mathbb{R}^{N}.
\end{align*}
Note that if we maintain the policy $\pi_{k}\triangleq\pi_{\theta_k}$ under the softmax parameterization, \eqref{Alg:AC:actor_update} is equivalent to \eqref{Alg:AC:actor_update_pi} with\looseness=-1
\begin{align*}
    \widehat{Z}_k(s)=\sum_{a' \in \Acal} \pi_k(a' \mid s) \exp \left(\alpha \widehat{Q}_{L,k}(s, a')\right)\in\mathbb{R}.
\end{align*}
In addition, the operator $\Pi$ in \eqref{Alg:AC:dual_update} defines an entry-wise projection to the interval $[0,\frac{2}{\xi(1-\gamma)}]$, motivated by the result in Lemma \ref{lem:bounded_lambdastar}.

The stochastic gradients of the primal and dual variables depend on the accuracy of the critic, which is updated according to \eqref{Alg:AC:critic_update}. This is an asynchronous Q table update using only one sample per iteration.

Essentially consisting of the critic update in \eqref{Alg:AC:critic_update}, the actor update in \eqref{Alg:AC:actor_update}, and the dual variable update in \eqref{Alg:AC:dual_update}, Algorithm~\ref{Alg:AC} is truly online in that it only needs one sample from the Markov chain controlled by $\widehat{\pi}_k$ in every iteration $k$, i.e., Eq.\ \eqref{Alg:AC:sample}. Defined in \eqref{Alg:AC:behaviorpolicy_update}, the behavior policy $\widehat{\pi}_k$ is the $\epsilon$-exploration version of the current actor policy $\pi_k$, which guarantees sufficient exploration in our algorithm (i.e., sufficient visitation to each state and action pair). Besides the step sizes $\eta,\alpha,\beta$ that are used for the dual, actor, and critic update, $\epsilon$ is another important parameter that needs to be carefully selected to guarantee the convergence of our algorithm. The same behavior policy was also considered in \citet{borkar2005actor,khodadadian2021finite}.


\begin{algorithm}[!h]
\SetAlgoLined
\textbf{Initialize:} $\pi_{0}$ and $\widehat{\pi}_0$ uniform distribution, $\widehat{Q}_{i,0} = 0\in\mathbb{R}^{|\Scal||\Acal|}$, and $\lambda_{0} = 0\in\Rset^{N}$. Draw the initial sample $s_0\sim\rho$ and $a_0\sim\widehat{\pi}_{0}(\cdot\mid s_0)$.




 \For{$k=0,1,\cdots,K$}{
    1) Observe $s_{k+1}\sim\Pcal(\cdot\mid s_k, a_k)$ and take action
    \begin{align}
    a_{k+1}\sim\widehat{\pi}_{k}(\cdot\mid s_{k+1}).\label{Alg:AC:sample}
    \end{align}
    
    2) Critic update: $\forall i=0,1,\cdots,N$
        \begin{align}
        &\widehat{Q}_{i,k+1}(s, a)=\widehat{Q}_{i,k}(s, a)+\label{Alg:AC:critic_update}\\
        &\hspace{5pt}\beta\hspace{-1pt}\big(r_i(s_{k}, \hspace{-1pt}a_{k})\hspace{-2pt}+\hspace{-2pt}\gamma \widehat{Q}_{i,k}(s_{k+1}, \hspace{-1pt}a_{k+1})\hspace{-2pt}-\hspace{-2pt}\widehat{Q}_{i,k}(s_k, \hspace{-1pt}a_{k})\hspace{-1pt}\big)\hspace{-1pt}.\notag
        \end{align}

    3) Actor update: $\forall s\in\Scal, a\in\Acal$
        \begin{align}
        \begin{aligned}
        \pi_{k+1}(a \mid s)&=\pi_k(a \mid s)\\
        &\hspace{15pt}\times\frac{\exp \left(\alpha \widehat{Q}_{L,k}(s, a)\right)}{\widehat{Z}_k(s)}.
        \end{aligned}
        \label{Alg:AC:actor_update_pi}
        \end{align}
    
    4) Behavior policy update: $\forall s\in\Scal,a\in\Acal$
    \begin{align}
        \widehat{\pi}_{k+1}(a\hspace{-2pt}\mid\hspace{-2pt} s)=\frac{\epsilon}{|\Acal|}+\left(1-\epsilon\right) \pi_{k+1}(a\hspace{-2pt}\mid\hspace{-2pt} s).
        \label{Alg:AC:behaviorpolicy_update}
    \end{align}

    5) Dual variable update:
    \begin{align}
        &\lambda_{k+1}=\label{Alg:AC:dual_update}\\
        &\hspace{15pt}\Pi\hspace{-1pt}\left(\hspace{-2pt}\lambda_k\hspace{-2pt}-\hspace{-2pt}\eta_k\Big(\sum_{s,a}\rho(s)\pi_k(a\hspace{-2pt}\mid\hspace{-2pt} s)\widehat{Q}_{g,k}(s,a)\hspace{-2pt}-\hspace{-2pt}b\Big)\hspace{-2pt}\right)\hspace{-2pt}.\notag
    \end{align}
 }
\caption{Online Primal-Dual Natural Actor-Critic Algorithm}
\label{Alg:AC}
\end{algorithm}

\section{Main Results}\label{sec:mainresults}

This section presents our finite-time analysis of Algorithm~\ref{Alg:AC} and highlights the main technical difficulty. In Section~\ref{sec:mainresults:convergence_rate}, we state our main theorem, which characterizes the convergence of the actor to the global optimality of \eqref{eq:objective} in both the objective function and constraint violation, under challenges caused by the coupled multi-time-scale updates and the time-varying Markovian samples. In Section~\ref{sec:mainresults:challenges}, we discuss these challenges in details, sketch our approach to treating the challenges, and introduce preliminary results on the convergence of the dual variable, actor, and critic as coupled functions of each other. Finally, we provide the proof of the main theorem in Section~\ref{sec:mainresults:proof_theorem}.

\subsection{Finite-Time Convergence}\label{sec:mainresults:convergence_rate}

We start by defining the mixing time of a Markov chain $\{X_k\}$, which measures the time for $\{X_k\}$ to approach its stationary distribution \citep{LevinPeresWilmer2006}. 


\begin{definition}\label{def:mixing_time}
Let $\mu$ denote the stationary distribution of the Markov chain $\{X_k\}$.
For any scalar $c>0$, the mixing time of $\{X_k\}$ associated with $c$ is defined as
\begin{align*}
    \tau(c) &=\\ &\hspace{-15pt}\min\{k\hspace{-2pt}\geq0:\hspace{-2pt}\sup_{X\in\Xcal}\hspace{-2pt}d_{\text{TV}}(P(X_k=\cdot\,|X_0\hspace{-2pt}=X),\mu(\cdot))\leq c\},
\end{align*}
where given two probability distributions $u_1$ and $u_2$, $d_{TV}$ denotes the total variation distance between them
\begin{equation}
    d_{\text{TV}}(u_1,u_2)=\frac{1}{2} \sup _{\nu: \Xcal \rightarrow[-1,1]}\left|\int \nu d u_1-\int \nu d u_2\right|.
    \label{eq:TV_def}
\end{equation}
\end{definition}
We consider the following assumption on the Markov chain generated by the transition probability $P^{\pi}(s' \mid s)=\sum_{a\in\Acal}\Pcal(s'\mid s,a)\pi(a\mid s)$ under the policy $\pi$, which basically says that it mixes geometrically fast. We use $\mu_{\pi}$ to denote the stationary distribution of the Markov chain under the transition probability $P^{\pi}$.
\begin{assump}\label{assump:markov-chain}
Given any $\pi$, the Markov chain $\{X_{k}\}$ generated by $P^{\pi}$ has a unique stationary distribution $\mu_{\pi}$, and is uniformly geometrically ergodic, i.e., there exist $C_{0}\geq 1$ and $\ell\in (0,1)$ such that
\begin{equation*}
\sup_{X\in\Xcal}d_{\text{TV}}(P^{\pi}(X_k=\cdot \mid X_0=X),\mu_{\pi}(\cdot))\hspace{-2pt}\leq\hspace{-2pt} C_{0}\ell^k,\, \forall k\geq 0.
\end{equation*}
\end{assump}
This assumption is standard in the analysis of RL methods under time-varying Markov chains \citep{zou2019finite,wu2020finite,khodadadian2021finite}. It also implies that there exists a positive constant $D$ such that
\begin{equation}
\tau(c) \leq D\log\left(1/c\right)\quad\forall c.\label{eq:mixing:tau}
\end{equation} 
In this work, we denote $\tau\triangleq\tau(\alpha)$, where $\alpha$ is the step size for the actor update. Due to the uniform ergodicity assumption, the stationary distribution $\mu_{\pi}$ is uniformly bounded away from 0 for any policy $\pi$, and we define
$\underline{\mu}\triangleq\min_{\pi,s}\mu_{\pi}(s)>0$.

\begin{thm}\label{thm:main}
Suppose that Assumptions~\ref{assump:Slater} and \ref{assump:markov-chain} hold. Consider the actor iterates $\{\pi_k\}$ obtained from $K$ iterations of Algorithm \ref{Alg:AC}. Let the step size sequences be
\begin{align}
    \eta=\frac{\eta_0}{K^{5/6}},\,\,\alpha=\frac{\alpha_0}{K^{5/6}},\,\, \beta=\frac{\beta_0}{K^{1/2}},\,\, \epsilon=\frac{\epsilon_0}{K^{1/6}},
    \label{thm:main:stepsize}
\end{align}
with $\frac{(1-\gamma)\underline{\mu}\epsilon_0\beta_0}{|\Acal|}\leq 1$. Then the objective function converges to the optimal value as 
\begin{align*}
    \frac{1}{K}\sum_{k=0}^{K-1}\mathbb{E}\left[V_{0}^{\pi^{\star}}(\rho)-V_{0}^{\pi_{k}}(\rho)\right]\leq\Ocal\big(\frac{N^2\log(K)}{K^{1/6}}\big),
\end{align*}
and the constraint violation decays to zero as
\begin{align*}
    \frac{1}{K}\sum_{k=0}^{K-1}\sum_{i=1}^{N}\mathbb{E}\left[\Big[b_i-V_{i}^{\pi_k}(\rho)\Big]_{+}\right]\leq \Ocal\big(\frac{N^2\log(K)}{K^{1/6}}\big).
\end{align*}
\end{thm}
In Theorem \ref{thm:main}, we show that Algorithm~\ref{Alg:AC} converges at a rate $\Ocal(1/K^{1/6})$. As we only draw one sample in every iteration of Algorithm~\ref{Alg:AC}, this theorem translates directly to a sample complexity of $\Ocal(1/\delta^6)$, i.e. at most $\Ocal(1/\delta^6)$ samples are needed for the objective function value and constraint violation to decay to $\delta$.
Compared with the natural policy gradient algorithm for CMDP with a $\Ocal(1/K^{1/2})$ convergence rate \citep{ding2020natural}, Algorithm~\ref{Alg:AC} needs to estimate the value functions online using additional time scales, which slows down its convergence. 
It has been observed in unconstrained MDPs that the convergence of two-time-scale actor-critic algorithms is also slower than their single-time-scale counterparts where an exact gradient estimate is assumed to be given \citep{wu2020finite,khodadadian2021finite}. Finally, we note that \citet{borkar2005actor} is the first to consider this online primal-dual actor-critic algorithm, where an asymptotic convergence is studied. Theorem \ref{thm:main} of this paper is the first known result on its finite-time convergence.

\subsection{Technical Challenges \& Solution Sketch}\label{sec:mainresults:challenges}

In this section, we discuss the major challenges in the proof of Theorem \ref{thm:main}, sketch our approach to the challenges, and introduce intermediate convergence results that will be later used in the proof of the theorem in Section \ref{sec:mainresults:proof_theorem}.

\noindent\textbf{Technical Challenges.} The first challenge in our analysis is to handle the multi-time-scale nature of the updates. Recall that Algorithm~\ref{Alg:AC} uses four step size parameters, namely, the actor step size $\alpha$, the critic step size $\beta$, the dual variable step size $\eta$, and the exploration parameter $\epsilon$. \citet{borkar2005actor} shows that the relative order $\eta\leq\alpha\leq\beta\leq\epsilon$ leads to the asymptotic convergence of the algorithm, but the exact choice of the parameters as a function of the horizon $K$ that yields a finite-time convergence rate remains unknown. In this work, we fill in this gap by properly selecting and balancing the step sizes through a four-time-scale analysis. Compared with the existing works in two-time-scale stochastic approximation/optimization \citep{doan2020nonlinear,hong2020two,zeng2021two,doan2021finite,Doan2019,chen2019finite} and analysis of actor-critic algorithms involving up to three time scales \citep{wu2020finite,khodadadian2021finite}, the additional time scale(s) of this work further complicates the analysis.

The second technical challenge of the work stems from the time-varying Markov samples. To deal with the noise of the dependent samples from time-invariant Markov processes, one can utilize the fast mixing time imposed by Assumption~\ref{assump:markov-chain} to show that the noise decays geometrically fast \citep{srikant2019finite,GuptaSY2019_twoscale,doan2020finite,zeng2020finite}. However, such analysis is not directly applicable here as the samples of our algorithm are from a time-varying Markov chain generated by the ever-changing behavior policy as illustrated below
\begin{align*}
    (s_0,a_0) \stackrel{\Pcal,\widehat{\pi}_{0}}{\longrightarrow} (s_1,a_1) \stackrel{\Pcal,\widehat{\pi}_{1}}{\longrightarrow} (s_2,a_2) \stackrel{\Pcal,\widehat{\pi}_{2}}{\longrightarrow} \hspace{-2pt}\cdots\hspace{-3pt} \stackrel{\Pcal,\widehat{\pi}_{k-1}}{\longrightarrow} (s_{k},a_{k}),
\end{align*}
which implies that the stationary distribution of this Markov chain changes with $\widehat{\pi}_k$ over time, and the noise caused by samples from this chain requires special treatment. In our analysis, we will show that the shift of the stationary distribution is controlled by the shift of the behavior policy, which can be properly bounded with the correct step size choice.

In addition, we note that the time-varying Markov processes above also imply that the variables depend on each other across time iterations. The updates of the critic and dual variable depend on samples generated by the $\epsilon$-exploration version of previous actor policies $\{\widehat{\pi}_t\}_{t\leq k-1}$ and in turn affect the actor update in the next iteration. A similar interaction also exists between the dual variable and the critic. Combined with the effect of multiple time scales, the variable coupling further creates increased difficulty of characterizing the convergence of each variable individually.


The key idea in our approach to handle the intertwined system of variables is to choose a proper decoupling scheme according to the structure of the updates. We observe that only the update of the critic variable in \eqref{Alg:AC:critic_update} directly utilizes the samples. The effect of the samples on the actor and dual variable is indirect through the critic. This motivates us to separate the analysis of the actor and dual variable from that of the critic. Specifically, we aim to establish 1) the convergence of the actor and the dual variable up to errors that depend on the accuracy of the critic and 2) the convergence of the critic as a function of the distribution shift of the Markov samples, which can be further traced back to the dependence on the one-step drift of the actor. In the sequel, we present more details of our approach.

\noindent\textbf{Convergence of Actor and Dual:} Since the actor and dual variable updates do not explicitly depend on the samples, we can avoid dealing with the complications arising from the Markovian samples in this part of the analysis. For now, we only focus on the convergence of actor and the dual variable by deriving the convergence of the objective function and constraint violation up to an error that scales with the inaccuracy of the critic variable. 
\begin{prop}\label{prop:convergence_obj_constraintviolation}
Under Assumption~\ref{assump:Slater}, the iterates of the actor $\{\pi_k\}$ converges to attain the optimal objective function value
\begin{align*}
    \frac{1}{K}\sum_{k=0}^{K-1}\left(V_{0}^{\pi^{\star}}(\rho)-V_{0}^{\pi_{k}}(\rho)\right)\hspace{-2pt}\leq\hspace{-2pt} \frac{\log|\Acal|}{(1-\gamma)K\alpha}\hspace{-2pt}+\hspace{-2pt}\frac{3N}{\xi(1-\gamma)^3 K} + \frac{4N\eta}{(1-\gamma)^3}+\frac{C_1}{K}\sum_{k=0}^{K-1}\sum_{i=0}^{N}\left\|\widehat{Q}_{i,k}-Q_{i}^{\pi_k}\right\|.&
\end{align*}
Furthermore, $\{\pi_k\}$ also converges in constraint violation
\begin{align*}
    &\frac{1}{K}\sum_{k=0}^{K-1}\sum_{i=1}^{N}\Big[b_i-V_{i}^{\pi_k}(\rho)\Big]_{+}\notag\\
    &\leq \frac{\xi(1-\gamma)}{2}\Big(\frac{\log|\Acal|}{(1-\gamma)K\alpha}+\frac{3N}{\xi(1-\gamma)^3 K}+\frac{4N\eta}{(1-\gamma)^3}+\frac{18N}{K\eta\xi^2(1-\gamma)^2}+\frac{C_1}{K}\sum_{k=0}^{K-1}\sum_{i=0}^{N}\|\hat{Q}_{i,k}-Q_{i}^{\pi_k}\|\Big),
\end{align*}
where $C_1=\frac{8(2-\gamma)}{\xi(1-\gamma)^3}+\frac{6}{\xi(1-\gamma)}$.
\end{prop}

Due to the space limit, we delay the proof of the proposition to the appendix. If the estimates of the critic were exactly correct, i.e. $\|\hat{Q}_{i,k}-Q_{i}^{\pi_k}\|=0$, we would obtain $\Ocal(1/\sqrt{K})$ convergence for both the objective function and constraint violation by choosing $\eta=\Ocal(1/\sqrt{K})$ and $\alpha=\Ocal(1)$, which recovers the convergence rate of the natural policy gradient method in \citet{ding2020natural}. In this work, such step size choices are not possible due to the presence of the critic error and the restrictions on the step sizes.

\noindent\textbf{Convergence of Critic:} We now establish the convergence of the critic to the true value function of the behavior policy, i.e. defining 
\begin{align}
    z_{i,k}=\widehat{Q}_{i,k}-Q_i^{\widehat{\pi}_k}\in\mathbb{R}^{|\Scal||\Acal|},\label{eq:def_z}
\end{align}
we aim to bound $\frac{1}{K}\sum_{k=0}^{K-1}\|z_{i,k}\|^2$.

As we have discussed earlier in this section, the coupling between the variables creates a dependency of the critic convergence on the step-wise drifts of the actor. We remove such dependency by choosing the step sizes such that the critic updates on a much faster time scale than the actor does. This essentially makes the drifts of the actor insignificant from the eye of the critic and reduces the bound on $\frac{1}{K}\sum_{k=0}^{K-1}\|z_{i,k}\|^2$ to only depend on the ratio of the step sizes.

\begin{prop}\label{prop:conv_critic}
Let Assumption \ref{assump:markov-chain} hold and
\[\frac{(1-\gamma)\underline{\mu}\epsilon\beta}{|\Acal|}\leq 1,\]
then the iterates $\{\widehat{Q}_{i,k}\}$ and $\{\widehat{\pi}_{k}\}$ satisfy
\begin{align*}
    \frac{1}{K}\sum_{k=0}^{K-1}\mathbb{E}\left[\|z_{i,k}\|^2\right]\leq \frac{8|\Scal||\Acal|^2 \tau}{(1-\gamma)^3\underline{\mu}\epsilon\beta K}+\frac{C_5 N^2|\Acal|\beta\tau^2}{(1-\gamma)\underline{\mu}\epsilon}+\frac{9|\Scal|^3|\Acal|^5 N^2 \alpha^2}{\xi\underline{\mu}^2 (1-\gamma)^9 \epsilon^2 \beta^2}&,
\end{align*}
where the constants are
\begin{align*}
    &C_2\hspace{-3pt}=\hspace{-3pt}\frac{2\sqrt{2}(|\Scal|\Acal|)^{3/2}}{(1-\gamma)^3}\hspace{-2pt}+\hspace{-2pt}\frac{8|\Scal|^2|\Acal|^3}{(1-\gamma)^2}\hspace{-1pt}\big(\hspace{-2pt}\left\lceil\log _{\rho} m^{-1}\hspace{-1pt}\right\rceil\hspace{-2pt}+\hspace{-2pt}\frac{1}{1\hspace{-2pt}-\hspace{-2pt}\rho}\hspace{-2pt}+\hspace{-2pt}2\big)\hspace{-1pt},\\
    &C_3\hspace{-3pt}=\hspace{-3pt}1\hspace{-2pt}+\hspace{-2pt}\frac{9\sqrt{2|\Scal||\Acal|}}{1-\gamma},\quad C_4\hspace{-3pt}=\hspace{-3pt}\frac{4|\Scal|^{3/2}|\Acal|^{3/2}}{1-\gamma}(1\hspace{-2pt}+\hspace{-2pt}\frac{3|\Scal||\Acal|}{1-\gamma}),\\
    &C_5\hspace{-3pt}=\hspace{-3pt}\frac{2(3C_2+6C_3+4C_4) (|\Scal||\Acal|)^{3/2} }{\xi^{1/2} (1-\gamma)^{7/2}}+\frac{36 (|\Scal||\Acal|)^3}{\xi(1-\gamma)^7}.
\end{align*}
\end{prop}

With its detailed proof presented in the appendix, this proposition relies on the fast mixing time of the Markov chain under all behavior policy iterates $\{\widehat{\pi}_k\}$ and the slow drift of the actor compared with the critic. We note that the critic learns the Q function of the behavior policy rather than that of the actor policy which we need to use in combination with Proposition~\ref{prop:convergence_obj_constraintviolation}, but the distance between the behavior policy and the actor policy can be properly controlled by the parameter $\epsilon$.

\subsection{Proof of Theorem~\ref{thm:main}}\label{sec:mainresults:proof_theorem}
In this section, we provide a proof of our main theorem using intermediate results established in Proposition~\ref{prop:convergence_obj_constraintviolation} and \ref{prop:conv_critic}. We first introduce the following lemma on the Lipschitz continuity of the Q function with respect to the policy.
\begin{lem}[\citet{khodadadian2021finite} Lemma 8]\label{lem:Q_Lipschitz_pi}
    For any policy $\pi_1$, $\pi_2$ and $i=0,1,\cdots,N$, we have
    \begin{align*}
        \|Q_i^{\pi_1}-Q_i^{\pi_2}\|\leq \frac{|\Scal||\Acal|}{(1-\gamma)^2}\|\pi_1-\pi_2\|.
    \end{align*}
\end{lem}

By Jensen's inequality, we have
\begin{align*}
\mathbb{E}[\|z_{i,k}\|]\leq(\mathbb{E}[\|z_{i,k}\|^2])^{1/2}.
\end{align*}
This implies
\begin{align}
    \frac{1}{K}\sum_{k=0}^{K-1}\mathbb{E}[\|z_{i,k}\|]&\leq \frac{1}{K}\sum_{k=0}^{K-1}(\mathbb{E}[\|z_{i,k}\|^2])^{1/2}\notag\\
    &\leq\left(\frac{1}{K}\sum_{k=0}^{K-1}\mathbb{E}[\|z_{i,k}\|^2]\right)^{1/2}\notag\\
    &\leq \Big(\frac{8|\Scal||\Acal|^2 \tau}{(1-\gamma)^3\underline{\mu}\epsilon\beta K}+\frac{C_5 N^2|\Acal|\beta\tau^2}{(1-\gamma)\underline{\mu}\epsilon}+\frac{9|\Scal|^3|\Acal|^5 N^2 \alpha^2}{\xi\underline{\mu}^2 (1-\gamma)^9 \epsilon^2 \beta^2}\Big)^{1/2}\notag\\
    &\leq \frac{2\sqrt{2}|\Scal|^{1/2}|\Acal|}{(1-\gamma)^{3/2}\underline{\mu}^{1/2}}\frac{\tau^{1/2}}{\epsilon^{1/2}\beta^{1/2} K^{1/2}}+\frac{C_5^{1/2} |\Acal|^{1/2}}{(1-\gamma)^{1/2}\underline{\mu}^{1/2}}\frac{N\beta^{1/2}\tau}{\epsilon^{1/2}}+\frac{3|\Scal|^{3/2}|\Acal|^{5/2}}{\xi^{1/2}\underline{\mu} (1-\gamma)^{9/2}}\frac{N\alpha}{\epsilon\beta}\notag\\
    &=\Ocal\left(\frac{\tau^{1/2}}{\epsilon^{1/2}\beta^{1/2} K^{1/2}}+\frac{N\beta^{1/2}\tau}{\epsilon^{1/2}}+\frac{N\alpha}{\epsilon\beta}\right),
    \label{thm:main:eq2}
\end{align}
where the second inequality follows from the Cauchy-Schwarz inequality, and the third inequality plugs in the bound from Proposition~\ref{prop:conv_critic}.

The bounds in Proposition \ref{prop:convergence_obj_constraintviolation} depend on the error $\frac{1}{K}\sum_{k=0}^{K-1}\sum_{i=0}^{N}\left\|\widehat{Q}_{i,k}-Q_{i}^{\pi_k}\right\|$, which we decompose and bound using \eqref{thm:main:eq2}. From the definition of $z_{i,k}$ in \eqref{eq:def_z}, we have for any $i=0,1,\cdots,N$
\begin{align}
    \frac{1}{K}\sum_{k=0}^{K-1}\mathbb{E}[\|\widehat{Q}_{i,k}-Q_{i}^{\pi_k}\|]&\leq \frac{1}{K}\sum_{k=0}^{K-1}\mathbb{E}[\|Q_{i}^{\widehat{\pi}_k}-Q_{i}^{\pi_k}\|]+\frac{1}{K}\sum_{k=0}^{K-1}\mathbb{E}[\|z_{i,k}\|]\notag\\
    &\leq \frac{1}{K}\sum_{k=0}^{K-1}\frac{|\Scal||\Acal|}{(1-\gamma)^2}\mathbb{E}[\|\widehat{\pi}_k-\pi_k\|]+\frac{1}{K}\sum_{k=0}^{K-1}\mathbb{E}[\|z_{i,k}\|]\notag\\
    &\leq \frac{2|\Scal|^{3/2}|\Acal| \epsilon}{(1-\gamma)^2}+\frac{1}{K}\sum_{k=0}^{K-1}\mathbb{E}[\|z_{i,k}\|]\notag\\
    &=\Ocal\left(\epsilon+\frac{\tau^{1/2}}{\epsilon^{1/2}\beta^{1/2} K^{1/2}}+\frac{N\beta^{1/2}\tau}{\epsilon^{1/2}}+\frac{N\alpha}{\epsilon\beta}\right),
    \label{thm:main:eq3}
\end{align}
where the second inequality uses Lemma \ref{lem:Q_Lipschitz_pi}, the last inequality is due to \eqref{thm:main:eq2}, and the third inequality follows from
\begin{align*}
    &\|\widehat{\pi}_k-\pi_k\|=\|\frac{\epsilon}{|\Acal|}\1_{|\Scal||\Acal|}+(1-\epsilon)\pi_k-\pi_k\|\notag\\
    &\leq\epsilon\|\frac{1}{|\Acal|}\1_{|\Scal||\Acal|}\|+\epsilon\|\pi_k\|\leq \epsilon\frac{|\Scal|^{1/2}}{|\Acal|^{1/2}}+\epsilon|\Scal|^{1/2}\hspace{-2pt}\leq\hspace{-2pt}2\epsilon|\Scal|^{1/2}.
\end{align*}

Using \eqref{thm:main:eq3} in Proposition \ref{prop:convergence_obj_constraintviolation}, we obtain the following bound on the optimality gap
\begin{align*}
    &\frac{1}{K}\sum_{k=0}^{K-1}\mathbb{E}\left[V_{0}^{\pi^{\star}}(\rho)-V_{0}^{\pi_{k}}(\rho)\right]\notag\\
    &\leq \frac{\log|\Acal|}{(1-\gamma)K\alpha}+\frac{3N}{\xi(1-\gamma)^3 K}+ \frac{4N\eta}{(1-\gamma)^3}+\frac{C_1}{K}\sum_{k=0}^{K-1}\sum_{i=0}^{N}\mathbb{E}\left[\left\|\widehat{Q}_{i,k}-Q_{i}^{\pi_k}\right\|\right] \notag\\
    &=\Ocal\Big(\frac{1}{K\alpha}+\frac{N}{K}+N\eta+N\epsilon+\frac{N\tau^{1/2}}{\hspace{-2pt}\epsilon^{1/2}\beta^{1/2} K^{1/2}}+\frac{N^2 \beta^{1/2}\tau}{\epsilon^{1/2}}+\frac{N^2 \alpha}{\epsilon\beta}\Big).
\end{align*}

For the constraint violation,
\begin{align*}
    &\frac{1}{K}\sum_{k=0}^{K-1}\sum_{i=1}^{N}\mathbb{E}\left[\Big[b_i-V_{i}^{\pi_k}(\rho)\Big]_{+}\right]\notag\\
    &\leq\frac{\xi(1-\gamma)}{2}\Big(\frac{\log|\Acal|}{(1-\gamma)K\alpha}+\frac{3N}{\xi(1-\gamma)^3 K}+\frac{4N\eta}{(1-\gamma)^3}+\frac{18N}{K\eta\xi^2(1-\gamma)^2}+\frac{C_1}{K}\sum_{k=0}^{K-1}\sum_{i=0}^{N}\mathbb{E}[\|\hat{Q}_{i,k}-Q_{i}^{\pi_k}\|]\Big)\notag\\
    &=\Ocal\big(\frac{1}{K\alpha}+\frac{N}{K}+N\eta+\frac{N}{K\eta}+N\epsilon+\frac{N\tau^{1/2}}{\hspace{-2pt}\epsilon^{1/2}\beta^{1/2} K^{1/2}}+\frac{N^2 \beta^{1/2}\tau}{\epsilon^{1/2}}+\frac{N^2 \alpha}{\epsilon\beta}\big).
\end{align*}

Plugging in the step sizes in \eqref{thm:main:stepsize} to the two inequalities above and recognizing from \eqref{eq:mixing:tau} that
\begin{align*}
    \tau^{1/2} \leq \tau\leq D\log(1/\alpha)=D\log(\frac{K^{5/6}}{\alpha_0})=\Ocal(\log(K))
\end{align*}
completes the proof.

\qed

\section{Numerical Simulations}

\begin{figure}[ht]
  \centering
  \includegraphics[width=.5\linewidth]{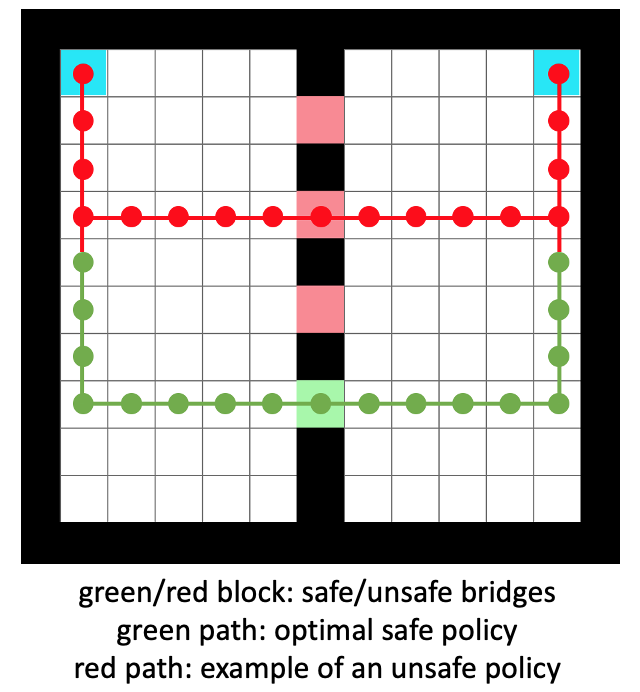}
  \caption{GridWorld}
  \label{fig:GridWorld}
\end{figure}

In this section, we experimentally validate the convergence of the online actor-critic algorithm on a GridWorld problem similar to the one considered in \citet{paternain2019constrained}, in which the agent aims to navigate to a goal position while obeying safety constraints in the environment illustrated in Figure~\ref{fig:GridWorld}. Starting from the top left corner, the agent receives a reward $r_0$ of $10$ if it reaches the goal position on the top right corner and a reward of $-1$ for in any other states. The choice of this reward function encourages the agent to find the shortest path to the goal position. At the same time, ``unsafe'' paths are defined through three constraints: we choose the utility functions $r_1$/$r_2$/$r_3$ to return $-1$ when the agent crosses the first/second/third red bridge in the second/fourth/sixth row and $0$ everywhere else, and constrain the values of the learned policy under the three utility functions to be larger than $-0.01$. The combined effect of the reward function and the constraints essentially requires the optimal policy to avoid the first three bridges and use the fourth one.

Applying Algorithm \ref{Alg:AC} to this problem, we observe that the learned policy indeed follows the optimal policy indicated by the green path in Figure~\ref{fig:GridWorld}. We numerically plot the converges of the empirical duality gap $\{g_k\}$ defined as
\begin{align*}
    g_k \triangleq \max_{\pi}V_L^{\pi,\lambda_k}(\rho)-\min_{\lambda\in\mathbb{R}_{+}^N}V_L^{\pi_k,\lambda}(\rho).
\end{align*}
It is easy to verify that $g_k=\max_{\pi}V_L^{\pi,\lambda_k}(\rho)-V_0^{\pi_k}(0)$ if $\pi_k$ is feasible and $g_k=\infty$ if $\pi_k$ is not. We initialize the actor policy to a feasible policy to avoid making the initial gap infinitely large. We note that in this small-scale GridWorld problem, the quantities $\max_{\pi}V_L^{\pi,\lambda_k}(\rho)$ and $V_0^{\pi_k}(0)$ can be quickly and accurately computed using dynamic programming.

\begin{figure}[ht]
  \centering
  \includegraphics[width=.6\linewidth]{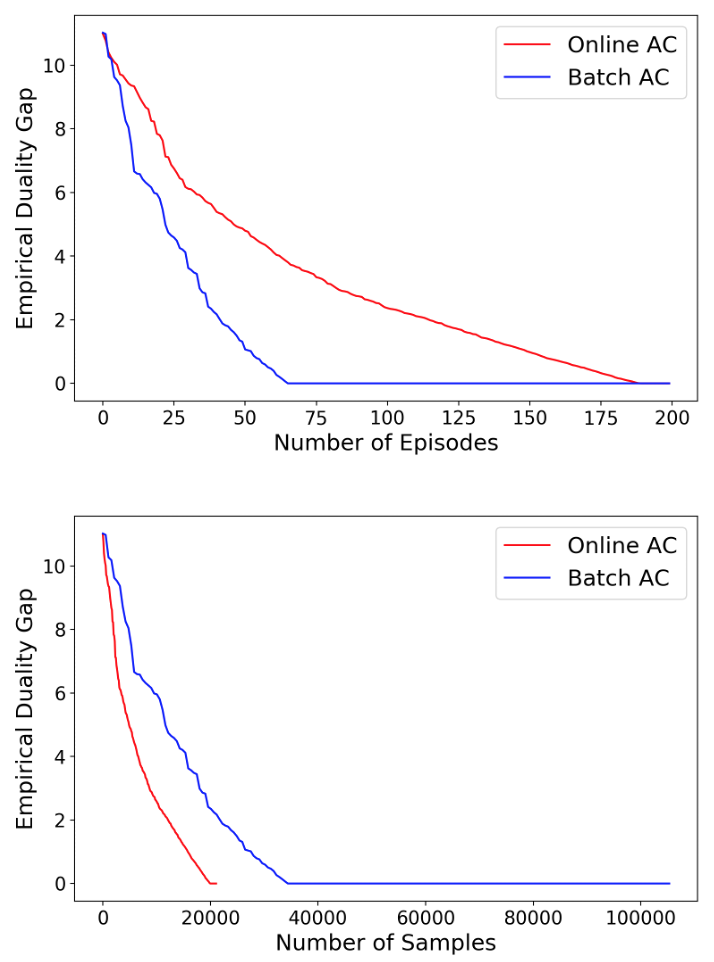}
  \caption{Duality Gap}
  \label{fig:duality_gap}
\end{figure}

We compare Algorithm \ref{Alg:AC} with a batch actor-critic algorithm where for every actor update, 5 separately generated trajectories are used to update the critic. We expect the batch algorithm to converge faster in the number of iterations but not to be as sample efficient, which is indeed observed in Figure~\ref{fig:duality_gap}.

\section{Conclusion \& Future Work}

In this paper, we derive the first-known finite-time global convergence of the online actor-critic algorithm for solving a CMDP problem. Future directions of the work include investigating the possibility of improving the convergence rate and extending the analysis to the scenario where function approximations are used to represent the actor and/or critic.

\newpage

\medskip
\bibliography{references}

\newpage
\onecolumn

    

\appendix

\section{Proof of Main Results}

We introduce a few shorthand notations that will be used in the rest of the paper. First, we define $\widehat{V}_{i,k}\in\mathbb{R}^{|\Scal|}$ and $\widehat{A}_{i,k}\in\mathbb{R}^{|\Scal|\times|\Acal|}$ such that
\begin{align*}
    \widehat{V}_{i,k}(s)=\sum_{a'\in\Acal}\pi_k(a'\mid s)\widehat{Q}_{i,k}(s,a'),\quad\text{and}\quad \widehat{A}_{i,k}(s,a) = \widehat{Q}_{i,k}(s,a)-\widehat{V}_{i,k}(s),\quad\forall s\in\Scal,a\in\Acal.
\end{align*}

We use the subscript $g$ under a value/Q/advantage function to denote the vector formed by stacking the value/Q/advantage functions over the utilities, i.e.
\begin{align*}
    &V_{g}^{\pi}(s) \triangleq [V_{1}^{\pi}(s),...,V_{N}^{\pi}(s)]^{\top}\in\mathbb{R}^{N}\\
    &\widehat{V}_{g,k}(s) \triangleq [\widehat{V}_{1,k}(s),...,\widehat{V}_{N,k}(s)]^{\top}\in\mathbb{R}^{N}\\
    &Q_{g}^{\pi}(s,a) \triangleq [Q_{1}^{\pi}(s,a),...,Q_{N}^{\pi}(s,a)]^{\top}\in\mathbb{R}^{N}\\
    &\widehat{Q}_{g,k}(s,a) \triangleq [\widehat{Q}_{1,k}(s,a),...,\widehat{Q}_{N,k}(s,a)]^{\top}\in\mathbb{R}^{N}\\
    &A_{g}^{\pi}(s,a) \triangleq [A_{1}^{\pi}(s,a),...,A_{N}^{\pi}(s,a)]^{\top}\in\mathbb{R}^{N}\\
    &\widehat{A}_{g,k}(s,a) \triangleq [\widehat{A}_{1,k}(s,a),...,\widehat{A}_{N,k}(s,a)]^{\top}\in\mathbb{R}^{N}.
\end{align*}

The subscript $L,k$ denotes the sum over reward and utilities of the value/Q/advantage functions, weighted by the dual variable in the $k$-th iteration, i.e.
\begin{align*}
    &V_{L,k}(s)\triangleq V_{0}^{\pi_k}(s)+\lambda_k^{\top} V_{g}^{\pi_k}(s)\in\mathbb{R}\\
    &\widehat{V}_{L,k}(s)\triangleq \widehat{V}_{0,k}(s)+\lambda_k^{\top} \widehat{V}_{g,k}(s)\in\mathbb{R}\\
    &Q_{L,k}(s, a)\triangleq Q_{0}^{\pi_k}(s, a)+\lambda_k^{\top} Q_{g}^{\pi_k}(s, a)\in\mathbb{R}\\
    &\widehat{Q}_{L,k}(s, a)\triangleq \widehat{Q}_{0,k}(s, a)+\lambda_k^{\top} \widehat{Q}_{g,k}(s, a)\in\mathbb{R}\\
    &A_{L,k}(s, a)\triangleq A_{0}^{\pi_k}(s, a)+\lambda_k^{\top} A_{g}^{\pi_k}(s, a)\in\mathbb{R}\\
    &\widehat{A}_{L,k}(s, a)\triangleq \widehat{A}_{0,k}(s, a)+\lambda_k^{\top} \widehat{A}_{g,k}(s, a)\in\mathbb{R}.
\end{align*}

Given a policy $\pi$ and an initial state distribution $\zeta$, we use $d_{\zeta}^{\pi}$ to denote the corresponding discounted visitation distribution such that
\begin{align}
    d_{\zeta}^{\pi}(s) = (1-\gamma)\mathbb{E}_{s_0\sim\zeta}\left[\sum_{k=0}^{\infty}\gamma^k P^{\pi}(s_k=s\mid s_0)\right].
    \label{eq:def_discounted_visitation}
\end{align}

We use $D_{\text{KL}}(p||q)$ to denote the Kullback–Leibler divergence between two distributions $p$ and $q$
\begin{align*}
    D_{\text{KL}}(p || q)\triangleq\mathbb{E}_{x \sim p} \log \frac{p(x)}{q(x)}.
\end{align*}

There is a simple bound we will frequently use that follows the boundedness of the reward and utility functions
\begin{align}
    |\lambda_{i,k}-\lambda_{i,k+1}|\leq \frac{2\eta}{1-\gamma},\,\forall i=1,2,\cdot,N, \,\,k=0,1,\cdots,K-1.\label{eq:lambda_ik_bound}
\end{align}



\subsection{Proof of Proposition \ref{prop:convergence_obj_constraintviolation}}

We start by introducing the following technical lemmas.

\begin{lem}\label{lem:bounded_Lagrangian}
The dual variable and actor iterates generated by Algorithm \ref{Alg:AC} satisfy
\begin{align*}
    &\frac{1}{K}\sum_{k=0}^{K-1}\hspace{-2pt}\left(\hspace{-2pt}V_{0}^{\pi^{\star}}(\rho)-V_{0}^{\pi_{k}}(\rho)+\lambda_k^{\top}\left(V_{g}^{\pi^{\star}}(\rho)-V_{g}^{\pi_{k}}(\rho)\hspace{-2pt}\right)\hspace{-2pt}\right)\notag\\
    &\leq\frac{\log|\Acal|}{(1-\gamma)K\alpha}+\frac{3N}{\xi(1-\gamma)^3 K}+\frac{2N\eta}{(1-\gamma)^3}+\frac{4-2\gamma}{(1-\gamma)^2 K}\sum_{k=0}^{K-1}\|A_{L,k}^{\pi_k}-\widehat{A}_{L,k}\|.
\end{align*}
\end{lem}

\begin{lem}\label{lem:step_improvement}
The dual variable and actor iterates generated by Algorithm \ref{Alg:AC} satisfy
\begin{align*}
    &V_{0}^{\pi_{k+1}}(\zeta)-V_{0}^{\pi_{k}}(\zeta)+\lambda_k^{\top}\left(V_{g}^{\pi_{k+1}}(s)-V_{g}^{\pi_k}(s)\right)\notag\\
    &\geq\frac{1}{\alpha}\mathbb{E}_{s \sim \zeta}\left[\log\widehat{Z}_k(s)\right]-\mathbb{E}_{s \sim \zeta}\left[\widehat{V}_{L,k}(s)\right]-\frac{3-\gamma}{1-\gamma}\|A_{L,k}^{\pi_k}-\widehat{A}_{L,k}\|
\end{align*}
for any initial distribution $\zeta$.
\end{lem}

Lemma \ref{lem:bounded_Lagrangian} bounds the
the difference between the empirical Lagrangian and that of the optimal solution averaged over iterations, while Lemma \ref{lem:step_improvement} considers the single-iteration improvement of the empirical Lagrangian.

\begin{lem}\label{lem:constraint_violation}
Suppose Assumption \ref{assump:Slater} holds.
Let the constant $C$ obey $C\geq 2\|\lambda^{\star}\|_{\infty}$.
Then, given a policy $\bar{\pi}$, if there exists a constant $\delta>0$ such that
\begin{align*}
    V_{0}^{\pi^{\star}}(\rho)-V_{0}^{\bar{\pi}}(\rho)+C\sum_{i=1}^{N}[b_i-V_i^{\bar{\pi}}(\rho)]_{+}\leq\delta,
\end{align*}
then we have
\begin{align*}
    \sum_{i=1}^{N}[b_i-V_i^{\bar{\pi}}(\rho)]_{+}\leq \frac{2\delta}{C}.
\end{align*}
\end{lem}

Lemma \ref{lem:constraint_violation} provides a convenient bound on the constraint violation of a given policy.

With the aid of the lemmas, we first study the convergence of the objective function to the global optimality.

\noindent\textbf{Objective function convergence}

From the dual update \eqref{Alg:AC:dual_update}, we have
\begin{align}
    0\leq\|\lambda_K\|^2&=\sum_{k=0}^{K-1}\left(\|\lambda_{k+1}\|^2-\|\lambda_k\|^2\right)\notag\\
    &=\sum_{k=0}^{K-1}\left(\left\|\Pi\left(\lambda_k-\eta\left(\sum_{s,a}\rho(s)\pi_k(a\mid s)\widehat{Q}_{g,k}(s,a)-b\right)\right)\right\|^2-\|\lambda_k\|^2\right)\notag\\
    &\leq\sum_{k=0}^{K-1}\left(\left\|\lambda_k-\eta\left(\sum_{s,a}\rho(s)\pi_k(a\mid s)\widehat{Q}_{g,k}(s,a)-b\right)\right\|^2-\|\lambda_k\|^2\right)\notag\\
    &=-2\eta\sum_{k=0}^{K-1}\lambda_k^{\top}\left(\sum_{s,a}\rho(s)\pi_k(a\mid s)\widehat{Q}_{g,k}(s,a)-b\right)+\eta^2\sum_{k=0}^{K-1}\left\|\sum_{s,a}\rho(s)\pi_k(a\mid s)\widehat{Q}_{g,k}(s,a)-b\right\|^2\notag\\
    &=-2\eta\sum_{k=0}^{K-1}\lambda_k^{\top}\left(\sum_{s,a}\rho(s)\pi_k(a\mid s)Q_{g}^{\pi_k}(s,a)-b\right)-2\eta\sum_{k=0}^{K-1}\lambda_k^{\top}\left(\sum_{s,a}\rho(s)\pi_k(a\mid s)(\widehat{Q}_{g,k}(s,a)-Q_{g}^{\pi_k}(s,a)\right)\notag\\
    &\hspace{20pt}+\eta^2\sum_{k=0}^{K-1}\left\|\sum_{s,a}\rho(s)\pi_k(a\mid s)\widehat{Q}_{g,k}(s,a)-b\right\|^2\notag\\
    &=-2\eta\sum_{k=0}^{K-1}\lambda_k^{\top}\left(V_g^{\pi_k}(\rho)-b\right)-2\eta\sum_{k=0}^{K-1}\lambda_k^{\top}\left(\sum_{s,a}\rho(s)\pi_k(a\mid s)(\widehat{Q}_{g,k}(s,a)-Q_{g}^{\pi_k}(s,a)\right)\notag\\
    &\hspace{20pt}+\eta^2\sum_{k=0}^{K-1}\left\|\sum_{s,a}\rho(s)\pi_k(a\mid s)\widehat{Q}_{g,k}(s,a)-b\right\|^2.
    \label{prop:convergence_obj_constraintviolation:eq0}
    \end{align}

By the Cauchy-Schwarz inequality the second term above can be bounded as
\begin{align}
    -\sum_{k=0}^{K-1}\lambda_k^{\top}\left(\sum_{s,a}\rho(s)\pi_k(a\mid s)(\widehat{Q}_{g,k}(s,a)-Q_{g}^{\pi_k}(s,a)\right)&=-\sum_{k=0}^{K-1}\sum_{i=1}^{N}\lambda_{i,k}\left(\sum_{s,a}\rho(s)\pi_k(a\mid s)(\widehat{Q}_{i,k}(s,a)-Q_{i}^{\pi_k}(s,a)\right)\notag\\
    &\leq \sum_{k=0}^{K-1}\sum_{i=1}^{N}\lambda_{i,k}\left(\sum_{s,a}\rho(s)^2\pi_k(a\mid s)^2\right)^{1/2}\|\widehat{Q}_{g,k}-Q_{g}^{\pi_k}\|\notag\\
    &\leq\sum_{k=0}^{K-1}\sum_{i=1}^{N}\lambda_{i,k}\left(\sum_{s,a}\rho(s)\pi_k(a\mid s)\right)\|\widehat{Q}_{g,k}-Q_{g}^{\pi_k}\|\notag\\
    &\leq\sum_{k=0}^{K-1}\sum_{i=1}^{N}\lambda_{i,k}\|\widehat{Q}_{g,k}-Q_{g}^{\pi_k}\|.
    \label{prop:convergence_obj_constraintviolation:eq0.25}
\end{align}
In addition, the third term of \eqref{prop:convergence_obj_constraintviolation:eq0} can be treated with the fact that the value function and constant $b_i$ are all within $[0,\frac{1}{1-\gamma}]$.
\begin{align}
    \sum_{k=0}^{K-1}\left\|\sum_{s,a}\rho(s)\pi_k(a\mid s)\widehat{Q}_{g,k}(s,a)-b\right\|^2&= \sum_{k=0}^{K-1}\sum_{i=1}^{N}\left\|\sum_{s,a}\rho(s)\pi_k(a\mid s)\widehat{Q}_{i,k}(s,a)-b_i\right\|^2\notag\\
    &\leq 2\sum_{k=0}^{K-1}\sum_{i=1}^{N}\left(\left\|\sum_{s,a}\rho(s)\pi_k(a\mid s)\widehat{Q}_{i,k}(s,a)\right\|^2+\left\|b_i\right\|^2\right)\notag\\
    &\leq 2\sum_{k=0}^{K-1}\sum_{i=1}^{N}\left(\frac{1}{(1-\gamma)^2}+\frac{1}{(1-\gamma)^2}\right)=\frac{4KN}{(1-\gamma)^2}.\label{prop:convergence_obj_constraintviolation:eq0.5}
\end{align}

Using \eqref{prop:convergence_obj_constraintviolation:eq0.25} and \eqref{prop:convergence_obj_constraintviolation:eq0.5} in \eqref{prop:convergence_obj_constraintviolation:eq0}, we get
\begin{align*}
    0&\leq-2\eta\sum_{k=0}^{K-1}\lambda_k^{\top}\left(V_{g}^{\pi_k}(\rho)-b\right)-2\eta\sum_{k=0}^{K-1}\lambda_k^{\top}\left(\sum_{s,a}\rho(s)\pi_k(a\mid s)(\widehat{Q}_{g,k}(s,a)-Q_{g}^{\pi_k}(s,a)\right)\notag\\
    &\hspace{20pt}+\eta^2\sum_{k=0}^{K-1}\left\|\sum_{s,a}\rho(s)\pi_k(a\mid s)\widehat{Q}_{g,k}(s,a)-b\right\|^2\notag\\
    &\leq-2\eta\sum_{k=0}^{K-1}\lambda_k^{\top}\left(V_{g}^{\pi_k}(\rho)-b\right)+2\eta\sum_{k=0}^{K-1}\sum_{i=1}^{N}\lambda_{i,k}\left\|\widehat{Q}_{i,k}-Q_{i}^{\pi_k}\right\|+\frac{4KN\eta^2}{(1-\gamma)^2}\notag\\
    &\leq 2\eta\sum_{k=0}^{K-1}\lambda_k^{\top} \left(V_{g}^{\pi^{\star}}(\rho)-V_{g}^{\pi_k}(\rho)\right)+2\eta\sum_{k=0}^{K-1}\sum_{i=1}^{N}\lambda_{i,k}\left\|\widehat{Q}_{i,k}-Q_{i}^{\pi_k}\right\|+\frac{4KN\eta^2}{(1-\gamma)^2},
\end{align*}
where the last inequality follows from the fact that the optimal policy satisfies the constraints, i.e. $V_i^{\pi^{\star}}(\rho)\geq b_i$ for all $i=1,\cdots,N$.

Re-arranging this inequality and dividing by $2K\eta$ lead to
\begin{align}
    \frac{1}{K}\sum_{k=0}^{K-1}\lambda_k^{\top} \left(V_{g}^{\pi^{\star}}(\rho)-V_{g}^{\pi_k}(\rho)\right)\geq -\frac{1}{K}\sum_{k=0}^{K-1}\sum_{i=1}^{N}\lambda_{i,k}\left\|\widehat{Q}_{i,k}-Q_{i}^{\pi_k}\right\| - \frac{2N\eta}{(1-\gamma)^2}.
\end{align}

Combining this with Lemma \ref{lem:bounded_Lagrangian}, we get
\begin{align}
    \frac{1}{K}\sum_{k=0}^{K-1}\left(V_{0}^{\pi^{\star}}(\rho)-V_{0}^{\pi_{k}}(\rho)\right)
    &\leq \frac{\log|\Acal|}{(1-\gamma)K\alpha}+\frac{3N}{\xi(1-\gamma)^3 K}+\frac{2N\eta}{(1-\gamma)^3}+\frac{4-2\gamma}{(1-\gamma)^2 K}\sum_{k=0}^{K-1}\|A_{L,k}^{\pi_k}-\widehat{A}_{L,k}\|\notag\\
    &\hspace{20pt}+\frac{1}{K}\sum_{k=0}^{K-1}\sum_{i=1}^{N}\lambda_{i,k}\left\|\widehat{Q}_{i,k}-Q_{i}^{\pi_k}\right\| + \frac{2N\eta}{(1-\gamma)^2}\notag\\
    &\leq \frac{\log|\Acal|}{(1-\gamma)K\alpha}+\frac{3N}{\xi(1-\gamma)^3 K}+\frac{4(2-\gamma)}{(1-\gamma)^2 K}\sum_{k=0}^{K-1}\|\widehat{Q}_{0,k}-Q_{0}^{\pi_k}\|\notag\\
    &\hspace{20pt}+\left(\frac{8(2-\gamma)}{\xi(1-\gamma)^3 K }+\frac{2}{\xi(1-\gamma)K}\right)\sum_{k=0}^{K-1}\sum_{i=1}^{N}\left\|\widehat{Q}_{i,k}-Q_{i}^{\pi_k}\right\| + \frac{4N\eta}{(1-\gamma)^3}
    \label{prop:convergence_obj_constraintviolation:eq1}
\end{align}
where the last inequality follows from the boundness of $\|\lambda_{i,k}\|$ due to the projection operator and
\begin{align}
    \sum_{k=0}^{K-1}\|A_{L,k}^{\pi_k}-\widehat{A}_{L,k}\|&\leq \sum_{k=0}^{K-1}\left(\|A_{0}^{\pi_k}-\widehat{A}_{0,k}\|+\sum_{i=1}^{N}|\lambda_{i,k}|\|A_{i}^{\pi_k}-\widehat{A}_{i,k}\|\right)\notag\\
    &\leq \sum_{k=0}^{K-1}\|A_{0}^{\pi_k}-\widehat{A}_{0,k}\|+\frac{2}{\xi(1-\gamma)}\sum_{k=0}^{K-1}\sum_{i=1}^{N}\|A_{i}^{\pi_k}-\widehat{A}_{i,k}\|\notag\\
    &\leq 2\sum_{k=0}^{K-1}\|Q_{0}^{\pi_k}-\widehat{Q}_{0,k}\|+\frac{4}{\xi(1-\gamma)}\sum_{k=0}^{K-1}\sum_{i=1}^{N}\|Q_{i}^{\pi_k}-\widehat{Q}_{i,k}\|.
    \label{eq:Q_bound_A}
\end{align}


Defining $C_1=\frac{8(2-\gamma)}{\xi(1-\gamma)^3}+\frac{4}{\xi(1-\gamma)}$ and combining terms in \eqref{prop:convergence_obj_constraintviolation:eq1}, we have
\begin{align*}
    \frac{1}{K}\sum_{k=0}^{K-1}\left(V_{0}^{\pi^{\star}}(\rho)-V_{0}^{\pi_{k}}(\rho)\right)
    &\leq \frac{\log|\Acal|}{(1-\gamma)K\alpha}+\frac{3N}{\xi(1-\gamma)^3 K}+ \frac{4N\eta}{(1-\gamma)^3}+\frac{C_1}{K}\sum_{k=0}^{K-1}\sum_{i=0}^{N}\left\|\widehat{Q}_{i,k}-Q_{i}^{\pi_k}\right\|.
\end{align*}

We now show the decay of the constraint violation to zero.

\textbf{Constraint violation convergence}

For any $\lambda\in[0,\frac{4}{(1-\gamma)\xi}]^{N}$, since the projection operator $\Pi$ is non-expansive, we have
\begin{align*}
    \|\lambda_{k+1}-\lambda\|^2 &= \|\Pi(\lambda_k-\eta(\hat{V}_{g,k}(\rho)-b))-\lambda\|^2\notag\\
    &\leq \|\lambda_k-\eta(\hat{V}_{g,k}(\rho)-b)-\lambda\|^2\notag\\
    &= \|\lambda_k-\lambda\|^2-2\eta(\lambda_k-\lambda)^{\top}(\hat{V}_{g,k}(\rho)-b)+\eta^2\|\hat{V}_{g,k}(\rho)-b\|^2\notag\\
    &\leq \|\lambda_k-\lambda\|^2-2\eta(\lambda_k-\lambda)^{\top}(\hat{V}_{g,k}(\rho)-V_{g}^{\pi_k}(\rho))-2\eta(\lambda_k-\lambda)^{\top}(V_{g}^{\pi_k}(\rho)-b)+\frac{4N\eta^2}{(1-\gamma)^2},
\end{align*}
where the last inequality bounds the quadratic term using an approach similar to \eqref{prop:convergence_obj_constraintviolation:eq0.5}.

Re-arranging the terms and summing up from $k=0$ to $K-1$, we get
\begin{align}
    \frac{2\eta}{K}\sum_{k=0}^{K-1}(\lambda_k-\lambda)^{\top}(V_{g}^{\pi_k}(\rho)-b)
    &\leq\frac{1}{K}\left(\|\lambda_0-\lambda\|^2-\|\lambda_K-\lambda\|^2\right)-\frac{2\eta}{K}\sum_{k=0}^{K-1}(\lambda_k-\lambda)^{\top}(\hat{V}_{g,k}(\rho)-V_{g}^{\pi_k}(\rho))+\frac{4N\eta^2}{(1-\gamma)^2}\notag\\
    &\leq \frac{1}{K}\|\lambda_0-\lambda\|^2+\frac{2\eta}{K}\sum_{k=0}^{K-1}\langle\lambda_k-\lambda,\hat{V}_{g,k}(\rho)-V_{g}^{\pi_k}(\rho)\rangle+\frac{4N\eta^2}{(1-\gamma)^2}.
    \label{eq:constraint_violation_1}
\end{align}

To bound the second term of this inequality, we consider
\begin{align}
    \frac{2\eta}{K}\sum_{k=0}^{K-1}\langle\lambda_k-\lambda,\hat{V}_{g,k}(\rho)-V_{g}^{\pi_k}(\rho)\rangle
    &\leq \left(\frac{2}{\xi(1-\gamma)}+\frac{4}{\xi(1-\gamma)}\right)\frac{2\eta}{K}\sum_{k=0}^{K-1}\|\hat{V}_{g,k}(\rho)-V_{g}^{\pi_k}(\rho)\|\notag\\
    &\leq \frac{12\eta}{\xi(1-\gamma)K}\sum_{k=0}^{K-1}\sum_{i=1}^{N}\|\hat{Q}_{i,k}-Q_{i}^{\pi_k}\|
    \label{eq:constraint_violation_2}
\end{align}

Using \eqref{eq:constraint_violation_2} in \eqref{eq:constraint_violation_1} and dividing both sides by $2\eta$, we get
\begin{align}
    \frac{1}{K}\sum_{k=0}^{K-1}(\lambda_k-\lambda)^{\top}(V_{g}^{\pi_k}(\rho)-b) &\leq \frac{1}{2K\eta}\|\lambda_0-\lambda\|^2+\frac{6}{\xi(1-\gamma)K}\sum_{k=0}^{K-1}\sum_{i=1}^{N}\|\hat{Q}_{i,k}-Q_{i}^{\pi_k}\|
    +\frac{2N\eta}{(1-\gamma)^2}
    \label{eq:constraint_violation_3}
\end{align}

Adding \eqref{eq:constraint_violation_3} to the result of Lemma \ref{lem:bounded_Lagrangian}, we have
\begin{align*}
    &\frac{1}{K}\sum_{k=0}^{K-1}\left(V_{0}^{\pi^{\star}}(\rho)-V_{0}^{\pi_{k}}(\rho)+\lambda_k^{\top}\left(V_{g}^{\pi^{\star}}(\rho)-b\right)+\lambda^{\top}(V_g^{\pi_k}(\rho)-b)\right)\notag\\
    &\leq \frac{\log|\Acal|}{(1-\gamma)K\alpha}+\frac{3N}{\xi(1-\gamma)^3 K}+\frac{2N\eta}{(1-\gamma)^3}+\frac{4-2\gamma}{(1-\gamma)^2 K}\sum_{k=0}^{K-1}\|A_{L,k}^{\pi_k}-\widehat{A}_{L,k}\|\notag\\
    &\hspace{20pt}+\frac{1}{2K\eta}\|\lambda_0-\lambda\|^2+\frac{6}{\xi(1-\gamma)K}\sum_{k=0}^{K-1}\sum_{i=1}^{N}\|\hat{Q}_{i,k}-Q_{i}^{\pi_k}\|
    +\frac{2N\eta}{(1-\gamma)^2}
\end{align*}

Combining the terms and using the fact that $V_g^{\pi^{\star}}(\rho)-b\geq 0$, this inequality implies
\begin{align}
    &\frac{1}{K}\sum_{k=0}^{K-1}\left(V_{0}^{\pi^{\star}}(\rho)-V_{0}^{\pi_{k}}(\rho)+\lambda^{\top}(b-V_g^{\pi_k}(\rho))\right)\notag\\
    &\leq \frac{\log|\Acal|}{(1-\gamma)K\alpha}+\frac{3N}{\xi(1-\gamma)^3 K}+\frac{4N\eta}{(1-\gamma)^3}+\frac{1}{2K\eta}\|\lambda_0-\lambda\|^2\notag\\
    &\hspace{20pt}+\frac{4-2\gamma}{(1-\gamma)^2 K}\sum_{k=0}^{K-1}\|A_{L,k}^{\pi_k}-\widehat{A}_{L,k}\|+\frac{6}{\xi(1-\gamma)K}\sum_{k=0}^{K-1}\sum_{i=1}^{N}\|\hat{Q}_{i,k}-Q_{i}^{\pi_k}\|\notag\\
    &\leq \frac{\log|\Acal|}{(1-\gamma)K\alpha}+\frac{3N}{\xi(1-\gamma)^3 K}+\frac{4N\eta}{(1-\gamma)^3}+\frac{1}{2K\eta}\|\lambda_0-\lambda\|^2\notag\\
    &\hspace{20pt}+\frac{4(2-\gamma)}{(1-\gamma)^2 K}\sum_{k=0}^{K-1}\|Q_{0}^{\pi_k}-\widehat{Q}_{0,k}\|+\left(\frac{8(2-\gamma)}{\xi(1-\gamma)^3 K}+\frac{6}{\xi(1-\gamma)K}\right)\sum_{k=0}^{K-1}\sum_{i=1}^{N}\|\hat{Q}_{i,k}-Q_{i}^{\pi_k}\|\notag\\
    &= \frac{\log|\Acal|}{(1-\gamma)K\alpha}+\frac{3N}{\xi(1-\gamma)^3 K}+\frac{4N\eta}{(1-\gamma)^3}+\frac{1}{2K\eta}\|\lambda_0-\lambda\|^2+\frac{C_1}{K}\sum_{k=0}^{K-1}\sum_{i=0}^{N}\left\|\widehat{Q}_{i,k}-Q_{i}^{\pi_k}\right\|,
    \label{eq:constraint_violation_4}
\end{align}
where the second inequality again follows from \eqref{eq:Q_bound_A}, and the equality uses the constant definition $C_1=\frac{8(2-\gamma)}{\xi(1-\gamma)^3}+\frac{6}{\xi(1-\gamma)}$.

Now, choosing $\lambda$ such that 
\begin{align*}
    \lambda_i = \begin{cases}\frac{4}{(1-\gamma)\xi}, & \text { if } b_i-V_i^{\pi_k}(\rho)\geq 0 \\ 0, & \text { else }\end{cases}
\end{align*}

Then, \eqref{eq:constraint_violation_4} leads to
\begin{align}
    &\frac{1}{K}\sum_{k=0}^{K-1}\left(V_{0}^{\pi^{\star}}(\rho)-V_{0}^{\pi_{k}}(\rho)\right)+\frac{1}{K}\sum_{k=0}^{K-1}\sum_{i=1}^{N}\frac{4}{(1-\gamma)\xi}\Big[b_i-V_{i}^{\pi_k}(\rho)\Big]_{+}\notag\\
    &\leq \frac{\log|\Acal|}{(1-\gamma)K\alpha}+\frac{3N}{\xi(1-\gamma)^3 K}+\frac{4N\eta}{(1-\gamma)^3}+\frac{1}{2K\eta}\cdot\frac{36N}{\xi^2(1-\gamma)^2}+\frac{C_1}{K}\sum_{k=0}^{K-1}\sum_{i=0}^{N}\left\|\widehat{Q}_{i,k}-Q_{i}^{\pi_k}\right\|
    \label{eq:constraint_violation_5}
\end{align}

Note that there always exists a policy $\bar{\pi}_K$ such that $d_{\rho}^{\bar{\pi}_K}=\frac{1}{K}\sum_{k=0}^{K-1}d_{\rho}^{\pi_k}$, which implies
\begin{align*}
    V_i^{\bar{\pi}_K}=\frac{1}{K}\sum_{k=0}^{K-1} V_i^{\pi_k}\quad\forall i=0,1,\cdots,N.
\end{align*}

As a result, \eqref{eq:constraint_violation_5} becomes
\begin{align}
    &\left(V_{0}^{\pi^{\star}}(\rho)-V_{0}^{\bar{\pi}_{K}}(\rho)\right)+\frac{4}{(1-\gamma)\xi}\sum_{i=1}^{N}\Big[b_i-V_{i}^{\bar{\pi}_{K}}(\rho)\Big]_{+}\notag\\
    &\leq \frac{\log|\Acal|}{(1-\gamma)K\alpha}+\frac{3N}{\xi(1-\gamma)^3 K}+\frac{4N\eta}{(1-\gamma)^3}+\frac{18N}{K\eta\xi^2(1-\gamma)^2}+\frac{C_1}{K}\sum_{k=0}^{K-1}\sum_{i=0}^{N}\left\|\widehat{Q}_{i,k}-Q_{i}^{\pi_k}\right\|
\end{align}

Recall that Lemma \ref{lem:bounded_lambdastar} states that $2\|\lambda^{\star}\|_{\infty}\leq\frac{4}{(1-\gamma)\xi}$. Applying Lemma \ref{lem:constraint_violation} with $C=\frac{4}{(1-\gamma)\xi}$, we have
\begin{align*}
    &\sum_{i=1}^{N}\Big[b_i-V_{i}^{\bar{\pi}_{K}}(\rho)\Big]_{+}\hspace{-2pt}\notag\\
    &\leq\hspace{-2pt}\frac{\xi(1-\gamma)}{2}\Big(\frac{\log|\Acal|}{(1-\gamma)K\alpha}+\frac{3N}{\xi(1-\gamma)^3 K}+\frac{4N\eta}{(1-\gamma)^3}+\frac{18N}{K\eta\xi^2(1-\gamma)^2}+\frac{C_1}{K}\sum_{k=0}^{K-1}\sum_{i=0}^{N}\|\hat{Q}_{i,k}-Q_{i}^{\pi_k}\|\Big).
\end{align*}

\qed

\subsection{Proof of Proposition \ref{prop:conv_critic}}

To begin with, we define some convenient shorthand notations used in the proof of Proposition \ref{prop:conv_critic}. First, we use $O_k$ to denote the data observation used for variable updates in iteration $k$, i.e.
\[O_k=(s_k,a_k,s_{k+1},a_{k+1}).\]
Second, we use $R_i:|\Scal|\times|\Acal|\times|\Scal|\times|\Acal|\rightarrow\mathbb{R}^{|\Scal||\Acal|}$ to denote the vector reward operator such that
\begin{align}
    [R_i(s,a,s',a')]_{m}= \begin{cases} r_i(s,a), & \text { if } m=(s,a), \\ 0, &\text{ otherwise.}\end{cases}
    \label{eq:def_vectorR}
\end{align}

We also define the operator $M:|\Scal|\times|\Acal|\times|\Scal|\times|\Acal|\rightarrow\mathbb{R}^{|\Scal||\Acal|\times|\Scal||\Acal|}$ such that
\begin{align}
    &\text{if $(s,a=s',a')$:}\quad
    \begin{aligned}
    [M(s,a,s',a')]_{m,n}=\begin{cases} \gamma-1, & \text{if } m=n=(s,a) \\ 0, & \text { otherwise}\end{cases}\end{aligned}\notag\\
    &\text{if $(s,a\neq s',a')$:}\quad \begin{aligned} [M(s,a,s',a')]_{m,n}=\begin{cases} -1, & \text { if } m=n=(s,a) \\
    \gamma, & \text { if, } m=(s,a),n=(s',a') \\ 
    0, & \text { otherwise}\end{cases}
    \end{aligned}
    \label{eq:def_matrixM}
\end{align}

For any policy $\pi$, vector $z\in\mathbb{R}^{|\Scal||\Acal|}$, and observation tuple $O\in\Scal\times\Acal\times\Scal\times\Acal$, we define
\begin{align*}
    &\widebar{M}^{\pi}=\mathbb{E}_{s\sim\mu_{\pi},a\sim\pi(\cdot\mid s),s'\sim\Pcal(\cdot\mid s,a),a'\sim\pi(\cdot\mid s')}[M(s,a,s',a')],\\
    &\Gamma_i(\pi,z,O)=z^{\top}(R_i(O)+M(O)Q_i^{\pi})+z^{\top}(M(O)-\widebar{M}^{\pi})z
\end{align*}
$\widebar{M}^{\pi}$ is an important operator which we will show is negative definite under any behavior policy $\widehat{\pi}_k$, and the operator $\Gamma_i$ plays an important role in bounding the noise caused by the time-varying Markov chain.

We also introduce the following technical lemmas. 
\begin{lem}\label{lem:hatpi_onestepdrift}
The error between $\widehat{\pi}_{k}$ and $\widehat{\pi}_{k+1}$ is bounded for all $k\geq 0$
\begin{align*}
    \|\widehat{\pi}_{k+1}-\widehat{\pi}_{k}\|\leq \frac{3(|\Scal||\Acal|)^{1/2} N \alpha}{\xi^{1/2}(1-\gamma)^{3/2}}.
\end{align*}
\end{lem}

\begin{lem}\label{lem:z_onestepdrift}
We have for all $i=0,\cdots,N$ and $k\geq 0$
\begin{align*}
    \|z_{i,k+1}-z_{i,k}\|^2\leq \frac{36 (|\Scal||\Acal|)^3 N^2 \beta^2}{\xi(1-\gamma)^7},\quad\text{and}\quad
    \|z_{i,k+1}-z_{i,k}\|\leq \frac{6 (|\Scal||\Acal|)^{3/2} N \beta}{\xi^{1/2}(1-\gamma)^{7/2}}.
\end{align*}
\end{lem}

Lemma \ref{lem:hatpi_onestepdrift} and \ref{lem:z_onestepdrift} bound the single-iteration drift of $\hat{\pi}_k$ and $z_{i,k}$, respectively.

\begin{lem}[\citet{khodadadian2021finite} Lemma 3]\label{lem:barM_negdef}
The matrix $\widebar{M}^{\widehat{\pi}_k}$ is negative definite for all $k\geq 0$
\begin{align*}
    x^{\top}\widebar{M}^{\widehat{\pi}_k}x\leq -\frac{(1-\gamma)\underline{\mu}\epsilon}{|\Acal|}\|x\|^2, \quad\forall x\in\mathbb{R}^{|\Scal||\Acal|},k\geq 0,
\end{align*}
where $\underline{\mu}=\min_{\pi,s}\mu_{\pi}(s)$ is a positive constant due to the uniform ergodicity of the Markov chain under any policy.
\end{lem}

As we have discussed earlier, the matrix $\widebar{M}^{\widehat{\pi}_k}$ is negative definition for all $k=0,1,\cdots,K-1$, which is formally established by Lemma \ref{lem:barM_negdef}.

\begin{lem}\label{lem:bound_Gamma}
Under Assumption \ref{assump:markov-chain}, we have for all $k\geq\tau$ and $i=0,\cdots,N$
\begin{align*}
    \mathbb{E}[\Gamma_i(\widehat{\pi}_k,z_{i,k},O_k)] \leq \frac{(3C_2+6C_3+4C_4) (|\Scal||\Acal|)^{3/2} N \beta\tau^2}{\xi^{1/2}(1-\gamma)^{7/2}},
\end{align*}
where $C_2=\frac{2\sqrt{2}|\Scal|^{3/2}|\Acal|^{3/2}}{(1-\gamma)^3}+\frac{8|\Scal|^2|\Acal|^3}{(1-\gamma)^2}\left(\left\lceil\log _{\rho} m^{-1}\right\rceil+\frac{1}{1-\rho}+2\right)$, $C_3=1+\frac{9\sqrt{2|\Scal||\Acal|}}{1-\gamma}$, and $C_4=\frac{4|\Scal|^{3/2}|\Acal|^{3/2}}{1-\gamma}(1+\frac{3|\Scal||\Acal|}{1-\gamma})$.
\end{lem}

Lemma \ref{lem:bound_Gamma} bounds an difficult and important cross term that later arises in the analysis due to the time-varying Markovian noise.

With the definition of $R_i$ and $M$ in \eqref{eq:def_vectorR} and \eqref{eq:def_matrixM}, the critic update \eqref{Alg:AC:critic_update} can be re-expressed in the vector form as
\begin{align*}
    \widehat{Q}_{i,k+1}=\widehat{Q}_{i,k}+\beta(R_i(O_k)+M(O_k)\widehat{Q}_{i,k}),
\end{align*}
which implies
\begin{align*}
    z_{i,k+1}-z_{i,k}=-(Q_i^{\widehat{\pi}_{k+1}}-Q_i^{\widehat{\pi}_{k}})+\beta(R_i(O_k)+M(O_k)Q_i^{\widehat{\pi}_{k}}+M(O_k)z_{i,k}).
\end{align*}

Re-writing $\|z_{i,k+1}\|^2-\|z_{i,k}\|^2$
\begin{align}
    \|z_{i,k+1}\|^2-\|z_{i,k}\|^2&=2z_{i,k}^{\top}(z_{i,k+1}-z_{i,k})+\|z_{i,k+1}-z_{i,k}\|^2\notag\\
    &=2z_{i,k}^{\top}(z_{i,k+1}-z_{i,k}-\beta\widebar{M}^{\widehat{\pi}_{k}}z_{i,k})+\|z_{i,k+1}-z_{i,k}\|^2+2\beta z_{i,k}^{\top}\widebar{M}^{\widehat{\pi}_{k}}z_{i,k}\notag\\
    &=2\beta\Gamma_i(\widehat{\pi}_k,z_{i,k},O_k)-2z_{i,k}^{\top}(Q_i^{\widehat{\pi}_{k+1}}-Q_i^{\widehat{\pi}_{k}})+\|z_{i,k+1}-z_{i,k}\|^2+2\beta z_{i,k}^{\top}\widebar{M}^{\widehat{\pi}_{k}}z_{i,k}\notag\\
    &\leq 2\beta\Gamma_i(\widehat{\pi}_k,z_{i,k},O_k)+2\|z_{i,k}\|\|Q_i^{\widehat{\pi}_{k+1}}-Q_i^{\widehat{\pi}_{k}}\|+\|z_{i,k+1}-z_{i,k}\|^2+2\beta z_{i,k}^{\top}\widebar{M}^{\widehat{\pi}_{k}}z_{i,k}.
    \label{prop:conv_critic:eq1}
\end{align}

Using Lemma \ref{lem:z_onestepdrift}-\ref{lem:bound_Gamma} to bound the first, third and fourth term of \eqref{prop:conv_critic:eq1} in expectation, we have for all $k\geq\tau$
\begin{align}
    &\mathbb{E}[\|z_{i,k+1}\|^2-\|z_{i,k}\|^2]\notag\\
    &\leq 2\beta\mathbb{E}[\Gamma_i(\widehat{\pi}_k,z_{i,k},O_k)]+2\mathbb{E}[\|z_{i,k}\|\|Q_i^{\widehat{\pi}_{k+1}}-Q_i^{\widehat{\pi}_{k}}\|]+\mathbb{E}[\|z_{i,k+1}-z_{i,k}\|^2]+2\beta \mathbb{E}[z_{i,k}^{\top}\widebar{M}^{\widehat{\pi}_{k}}z_{i,k}]\notag\\
    &\leq \frac{2(3C_2+6C_3+4C_4) (|\Scal||\Acal|)^{3/2} N \beta^2 \tau^2}{\xi^{1/2} (1-\gamma)^{7/2}} + 2\mathbb{E}[\|z_{i,k}\|\|Q_i^{\widehat{\pi}_{k+1}}-Q_i^{\widehat{\pi}_{k}}\|]\notag\\
    &\hspace{20pt}+\frac{36 (|\Scal||\Acal|)^3 N^2 \beta^2}{\xi(1-\gamma)^7}-\frac{2(1-\gamma)\underline{\mu}\epsilon\beta}{|\Acal|}\mathbb{E}[\|z_{i,k}\|^2].
    \label{prop:conv_critic:eq2}
\end{align}

The second term of \eqref{prop:conv_critic:eq2} can be simplified with Lemma \ref{lem:Q_Lipschitz_pi} and \ref{lem:hatpi_onestepdrift}
\begin{align*}
    2\mathbb{E}[\|z_{i,k}\|\|Q_i^{\widehat{\pi}_{k+1}}-Q_i^{\widehat{\pi}_{k}}\|] &\leq \frac{2|\Scal||\Acal|}{(1-\gamma)^2}\mathbb{E}[\|z_{i,k}\|\|\widehat{\pi}_{k+1}-\widehat{\pi}_{k}\|]\notag\\
    &\leq \frac{6(|\Scal||\Acal|)^{3/2} N \alpha}{\xi^{1/2}(1-\gamma)^{7/2}}\mathbb{E}[\|z_{i,k}\|]\notag\\
     &\leq
     \frac{(1-\gamma)\underline{\mu}\epsilon\beta}{|\Acal|}\mathbb{E}[\|z_{i,k}\|^2]+
     \left(\frac{3(|\Scal||\Acal|)^{3/2} N \alpha}{\xi^{1/2}(1-\gamma)^{7/2}}\right)^2\frac{|\Acal|}{(1-\gamma)\underline{\mu}\epsilon\beta}\notag\\
     &\leq \frac{(1-\gamma)\underline{\mu}\epsilon\beta}{|\Acal|}\mathbb{E}[\|z_{i,k}\|^2] +
     \frac{9|\Scal|^3|\Acal|^4 N^2 \alpha^2}{\xi\underline{\mu}(1-\gamma)^8\epsilon\beta},
\end{align*}
where the third inequality follows from the fact that $2c\leq \frac{c^2}{d}+d$ for any scalar $c,d>0$. Plugging this bound back to \eqref{prop:conv_critic:eq2}, we have for all $k\geq\tau$
\begin{align*}
    \mathbb{E}[\|z_{i,k+1}\|^2-\|z_{i,k}\|^2]&\leq \frac{2(3C_2+6C_3+4C_4) (|\Scal||\Acal|)^{3/2} N \beta^2 \tau^2}{\xi^{1/2} (1-\gamma)^{7/2}}+\frac{(1-\gamma)\underline{\mu}\epsilon\beta}{|\Acal|}\mathbb{E}[\|z_{i,k}\|^2]\notag\\
    &\hspace{20pt}+
    \frac{9|\Scal|^3|\Acal|^4 N^2 \alpha^2}{\xi\underline{\mu}(1-\gamma)^8\epsilon\beta}+\frac{36 (|\Scal||\Acal|)^3 N^2 \beta^2}{\xi(1-\gamma)^7}-\frac{2(1-\gamma)\underline{\mu}\epsilon\beta}{|\Acal|}\mathbb{E}[\|z_{i,k}\|^2]\notag\\
    &\leq -\frac{(1-\gamma)\underline{\mu}\epsilon\beta}{|\Acal|}\mathbb{E}[\|z_{i,k}\|^2] + C_5 N^2 \beta^2\tau^2 + \frac{9|\Scal|^3|\Acal|^4 N^2 \alpha^2}{\xi\underline{\mu}(1-\gamma)^8\epsilon\beta},
\end{align*}
where in the second inequality we re-combine the terms using $\tau\geq1$, $N\geq1$, and the definition of the constant $C_5=\frac{2(3C_2+6C_3+4C_4) (|\Scal||\Acal|)^{3/2} }{\xi^{1/2} (1-\gamma)^{7/2}}+\frac{36 (|\Scal||\Acal|)^3}{\xi(1-\gamma)^7}$.

Re-arranging this inequality, we have for all $k\geq\tau$
\begin{align}
    \frac{(1-\gamma)\underline{\mu}\epsilon\beta}{|\Acal|}\mathbb{E}[\|z_{i,k}\|^2]&\leq \mathbb{E}[\|z_{i,k}\|^2-\|z_{i,k+1}\|^2] + C_5 N^2 \beta^2\tau^2 + \frac{9|\Scal|^3|\Acal|^4 N^2 \alpha^2}{\xi\underline{\mu}(1-\gamma)^8\epsilon\beta}.
    \label{prop:conv_critic:eq3}
\end{align}

As a result,
\begin{align*}
    &\frac{1}{K}\sum_{k=0}^{K-1}\frac{(1-\gamma)\underline{\mu}\epsilon\beta}{|\Acal|}\mathbb{E}[\|z_{i,k}\|^2]\notag\\
    &=\frac{1}{K}\sum_{k=0}^{\tau-1}\frac{(1-\gamma)\underline{\mu}\epsilon\beta}{|\Acal|}\mathbb{E}[\|z_{i,k}\|^2]+\frac{1}{K}\sum_{k=\tau}^{K-1}\frac{(1-\gamma)\underline{\mu}\epsilon\beta}{|\Acal|}\mathbb{E}[\|z_{i,k}\|^2]\notag\\
    &\leq \frac{1}{K}\sum_{k=0}^{\tau-1}\mathbb{E}[\|z_{i,k}\|^2]+\frac{1}{K}\sum_{k=\tau}^{K-1}\mathbb{E}[\|z_{i,k}\|^2-\|z_{i,k+1}\|^2] + \frac{1}{K}\sum_{k=\tau}^{K-1}\left(C_5 N^2 \beta^2\tau^2 + \frac{9|\Scal|^3|\Acal|^4 N^2 \alpha^2}{\xi\underline{\mu}(1-\gamma)^8\epsilon\beta}\right)\notag\\
    &\leq\frac{1}{K}\sum_{k=0}^{\tau}\mathbb{E}[\|z_{i,k}\|^2]+C_5 N^2 \beta^2\tau^2 + \frac{9|\Scal|^3|\Acal|^4 N^2 \alpha^2}{\xi\underline{\mu}(1-\gamma)^8\epsilon\beta},
\end{align*}
where the first inequality follows from \eqref{prop:conv_critic:eq3} and the step size condition $\frac{(1-\gamma)\underline{\mu}\epsilon\beta}{|\Acal|}\leq 1$.
Since $\|z_{i,k}\|^2=\|\widehat{Q}_{i,k}-Q_i^{\widehat{\pi}_k}\|^2\leq \frac{4|\Scal||\Acal|}{(1-\gamma)^2}$ for all $k\geq 0$, we can further simplify this inequality
\begin{align*}
    \frac{1}{K}\sum_{k=0}^{K-1}\frac{(1-\gamma)\underline{\mu}\epsilon\beta}{|\Acal|}\mathbb{E}[\|z_{i,k}\|^2] & \leq\frac{1}{K}\sum_{k=0}^{\tau}\mathbb{E}[\|z_{i,k}\|^2] + C_5 N^2 \beta^2\tau^2 + \frac{9|\Scal|^3|\Acal|^4 N^2 \alpha^2}{\xi\underline{\mu}(1-\gamma)^8\epsilon\beta}\notag\\
    &\leq\frac{4|\Scal||\Acal|(\tau+1)}{(1-\gamma)^2 K}+C_5 N^2 \beta^2\tau^2 + \frac{9|\Scal|^3|\Acal|^4 N^2 \alpha^2}{\xi\underline{\mu}(1-\gamma)^8\epsilon\beta}.
\end{align*}

Dividing both sides of the inequality by $\frac{(1-\gamma)\underline{\mu}\epsilon\beta}{|\Acal|}$, we get
\begin{align*}
    \frac{1}{K}\sum_{k=0}^{K-1}\mathbb{E}[\|z_{i,k}\|^2]&\leq \frac{|\Acal|}{(1-\gamma)\underline{\mu}\epsilon\beta}\frac{8|\Scal||\Acal|\tau}{(1-\gamma)^2 K}+\frac{|\Acal|}{(1-\gamma)\underline{\mu}\epsilon\beta}\left(C_5 N^2 \beta^2\tau^2 + \frac{9|\Scal|^3|\Acal|^4 N^2 \alpha^2}{\xi\underline{\mu}(1-\gamma)^8\epsilon\beta}\right)\notag\\
    &\leq \frac{8|\Scal||\Acal|^2 \tau}{(1-\gamma)^3\underline{\mu}\epsilon\beta K}+\frac{C_5 N^2|\Acal|}{(1-\gamma)\underline{\mu}}\frac{\beta\tau^2}{\epsilon}+\frac{9|\Scal|^3|\Acal|^5 N^2 \alpha^2}{\xi\underline{\mu}^2 (1-\gamma)^9 \epsilon^2 \beta^2}.
\end{align*}

\qed

\section{Proof of Technical Lemmas}

\subsection{Proof of Lemma \ref{lem:bounded_lambdastar}}
Let $\bar{\pi}$ denote the policy that satisfies the Slater's condition in Assumption \ref{assump:Slater}.
By the definition of the dual function,
\begin{align}
    V_D^{\lambda^{\star}}(\rho)\geq V_0^{\bar{\pi}}+(V_g^{\bar{\pi}}-b)^{\top}\lambda^{\star}\geq V_0^{\bar{\pi}}+\xi\sum_{i=1}^{N}\lambda_i^{\star}\geq V_0^{\bar{\pi}}+\xi\|\lambda^{\star}\|_{\infty}
\end{align}

Re-arranging the inequality and recognizing $V_D^{\lambda^{\star}}(\rho)=V_0^{\pi^{\star}}$ as implied by the strong duality in \eqref{lem:strong_duality} leads to
\begin{align*}
    \|\lambda^\star\|_{\infty}\leq \frac{1}{\xi}(V_D^{\lambda^{\star}}(\rho)-V_0^{\bar{\pi}})=\frac{1}{\xi}(V_0^{\pi^{\star}}(\rho)-V_0^{\bar{\pi}}(\rho))\leq \frac{2}{\xi(1-\gamma)},
\end{align*}
where the last inequality follows from bounded reward function (and hence bounded value function).

\qed





\subsection{Proof of Lemma \ref{lem:bounded_Lagrangian}}

The performance difference lemma \citep{agarwal2020optimality} states that for any policies $\pi_1$, $\pi_2$, initial distribution $\zeta$, and $i=0,\cdots,N$
\begin{align}
    V_i^{\pi_1}(\zeta)-V_i^{\pi_2}(\zeta)=\frac{1}{1-\gamma}\mathbb{E}_{s \sim d_{\zeta}^{\pi^{\star}}, a \sim \pi^{\star}(\cdot \mid s)}\left[A_{0}^{\pi_k}(s, a)\right].
    \label{eq:def_lemma_performancediff}
\end{align}

Using this lemma, we have
\begin{align*}
V_{0}^{\pi^{\star}}(\rho)-V_0^{\pi_k}(\rho)&= \frac{1}{1-\gamma}\mathbb{E}_{s \sim d_{\rho}^{\pi^{\star}}, a \sim \pi^{\star}(\cdot \mid s)}\left[A_{0}^{\pi_k}(s, a)\right] \notag\\
&=\frac{1}{1-\gamma}\mathbb{E}_{s \sim d_{\rho}^{\pi^{\star}},a\sim\pi^{\star}(\cdot\mid s)}\left[ Q_{0}^{\pi_k}(s, a)\right]-\frac{1}{1-\gamma}\mathbb{E}_{s \sim d_{\rho}^{\pi^{\star}}}\left[V_{0}^{\pi_k}(s)\right]\notag\\
&=\frac{1}{1-\gamma}\mathbb{E}_{s \sim d_{\rho}^{\pi^{\star}}, a \sim \pi^{\star}(\cdot \mid s)}\left[\widehat{Q}_{L,k}(s, a)\right]-\frac{\lambda_k^{\top}}{1-\gamma}\mathbb{E}_{s \sim d_{\rho}^{\pi^{\star}}, a \sim \pi^{\star}(\cdot \mid s)}\left[\widehat{Q}_{g,k}(s, a)\right]\notag\\
&\hspace{20pt}+\frac{1}{1-\gamma}\mathbb{E}_{s \sim d_{\rho}^{\pi^{\star}}, a \sim \pi^{\star}(\cdot \mid s)}\left[Q_{0}^{\pi_k}(s,a)-\widehat{Q}_{0,k}(s, a)\right]-\frac{1}{1-\gamma}\mathbb{E}_{s \sim d_{\rho}^{\pi^{\star}}}\left[V_{0}^{\pi_k}(s)\right] \notag\\
&=\frac{1}{(1-\gamma)\alpha}\mathbb{E}_{s \sim d_{\rho}^{\pi^{\star}}, a \sim \pi^{\star}(\cdot \mid s)}\left[\log\left(\frac{\pi_{k+1}(a\mid s)}{\pi_k(a\mid s)}\widehat{Z}_k\right)\right]\notag\\
&\hspace{20pt}-\frac{\lambda_k^{\top}}{1-\gamma}\mathbb{E}_{s \sim d_{\rho}^{\pi^{\star}}, a \sim \pi^{\star}(\cdot \mid s)}\left[\widehat{A}_{g,k}(s, a)\right]-\frac{\lambda_k^{\top}}{1-\gamma}\mathbb{E}_{s \sim d_{\rho}^{\pi^{\star}}}\left[\widehat{V}_{g,k}(s)\right]\notag\\
&\hspace{20pt}+\frac{1}{1-\gamma}\mathbb{E}_{s \sim d_{\rho}^{\pi^{\star}}, a \sim \pi^{\star}(\cdot \mid s)}\left[Q_{0}^{\pi_k}(s,a)-\widehat{Q}_{0,k}(s, a)\right]-\frac{1}{1-\gamma}\mathbb{E}_{s \sim d_{\rho}^{\pi^{\star}}}\left[V_{0}^{\pi_k}(s)\right] \notag\\
&=\frac{1}{(1-\gamma)\alpha}\mathbb{E}_{s \sim d_{\rho}^{\pi^{\star}}}\left[D_{\text{KL}}(\pi^{\star}(\cdot\mid s)||\pi_k(\cdot\mid s))\right]-\frac{1}{(1-\gamma)\alpha}\mathbb{E}_{s \sim d_{\rho}^{\pi^{\star}}}\left[D_{\text{KL}}(\pi^{\star}(\cdot\mid s)||\pi_{k+1}(\cdot\mid s))\right]\notag\\
&\hspace{20pt}+\frac{1}{(1-\gamma)\alpha}\mathbb{E}_{s \sim d_{\rho}^{\pi^{\star}}}\left[\log\widehat{Z}_k\right]-\frac{\lambda_k^{\top}}{1-\gamma}\mathbb{E}_{s \sim d_{\rho}^{\pi^{\star}}, a \sim \pi^{\star}(\cdot \mid s)}\left[\widehat{A}_{g,k}(s, a)\right]\notag\\
&\hspace{20pt}-\frac{\lambda_k^{\top}}{1-\gamma}\mathbb{E}_{s \sim d_{\rho}^{\pi^{\star}}}\left[\widehat{V}_{g,k}(s)\right]+\frac{1}{1-\gamma}\mathbb{E}_{s \sim d_{\rho}^{\pi^{\star}}, a \sim \pi^{\star}(\cdot \mid s)}\left[Q_{0}^{\pi_k}(s,a)-\widehat{Q}_{0,k}(s, a)\right]\notag\\
&\hspace{20pt}-\frac{1}{1-\gamma}\mathbb{E}_{s \sim d_{\rho}^{\pi^{\star}}}\left[\widehat{V}_{0,k}(s)\right]-\frac{1}{1-\gamma}\mathbb{E}_{s \sim d_{\rho}^{\pi^{\star}}}\left[V_{0}^{\pi_k}(s)-\widehat{V}_{0,k}(s)\right].
\end{align*}

Re-grouping the terms,
\begin{align}
V_{0}^{\pi^{\star}}(\rho)-V_0^{\pi_k}(\rho)&=\frac{1}{(1-\gamma)\alpha}\mathbb{E}_{s \sim d_{\rho}^{\pi^{\star}}}\left[D_{\text{KL}}(\pi^{\star}(\cdot\mid s)||\pi_k(\cdot\mid s))\right]-\frac{1}{(1-\gamma)\alpha}\mathbb{E}_{s \sim d_{\rho}^{\pi^{\star}}}\left[D_{\text{KL}}(\pi^{\star}(\cdot\mid s)||\pi_{k+1}(\cdot\mid s))\right]\notag\\
&\hspace{20pt}+\frac{1}{(1-\gamma)\alpha}\mathbb{E}_{s \sim d_{\rho}^{\pi^{\star}}}\left[\log\widehat{Z}_k\right]-\frac{\lambda_k^{\top}}{1-\gamma}\mathbb{E}_{s \sim d_{\rho}^{\pi^{\star}}, a \sim \pi^{\star}(\cdot \mid s)}\left[A_{g}^{\pi_k}(s, a)\right]\notag\\
&\hspace{20pt}+\frac{\lambda_k^{\top}}{1-\gamma}\mathbb{E}_{s \sim d_{\rho}^{\pi^{\star}}, a \sim \pi^{\star}(\cdot \mid s)}\left[A_{g}^{\pi_k}(s, a)-\widehat{A}_{g,k}(s, a)\right]\notag\\
&\hspace{20pt}-\frac{1}{1-\gamma}\mathbb{E}_{s \sim d_{\rho}^{\pi^{\star}}}\left[\widehat{V}_{L,k}(s)\right]+\frac{1}{1-\gamma}\mathbb{E}_{s \sim d_{\rho}^{\pi^{\star}}, a \sim \pi^{\star}(\cdot \mid s)}\left[A_{0}^{\pi_k}(s,a)-\widehat{A}_{0,k}(s, a)\right]\notag\\
&=\frac{1}{(1-\gamma)\alpha}\mathbb{E}_{s \sim d_{\rho}^{\pi^{\star}}}\left[D_{\text{KL}}(\pi^{\star}(\cdot\mid s)||\pi_k(\cdot\mid s))\right]-\frac{1}{(1-\gamma)\alpha}\mathbb{E}_{s \sim d_{\rho}^{\pi^{\star}}}\left[D_{\text{KL}}(\pi^{\star}(\cdot\mid s)||\pi_{k+1}(\cdot\mid s))\right]\notag\\
&\hspace{20pt}+\frac{1}{(1-\gamma)\alpha}\mathbb{E}_{s \sim d_{\rho}^{\pi^{\star}}}\left[\log\widehat{Z}_k\right]-\lambda_k^{\top}\left(V_{g}^{\pi^{\star}}(\rho)-V_{g}^{\pi_k}(\rho)\right)\notag\\
&\hspace{20pt}-\frac{1}{1-\gamma}\mathbb{E}_{s \sim d_{\rho}^{\pi^{\star}}}\left[\widehat{V}_{L,k}(s)\right]+\frac{1}{1-\gamma}\mathbb{E}_{s \sim d_{\rho}^{\pi^{\star}}, a \sim \pi^{\star}(\cdot \mid s)}\left[A_{L,k}^{\pi_k}(s, a)-\widehat{A}_{L,k}(s, a)\right]
\label{lem:bounded_Lagrangian:eq1}
\end{align}

Using Lemma \ref{lem:step_improvement} with $\zeta=d_{\rho}^{\pi^{\star}}$, we get
\begin{align}
    &\frac{1}{\alpha}\mathbb{E}_{s \sim d_{\rho}^{\pi^{\star}}}\left[\log\widehat{Z}_k(s)\right]-\mathbb{E}_{s \sim d_{\rho}^{\pi^{\star}}}\left[\widehat{V}_{L,k}(s)\right]\notag\\
    &\leq V_{0}^{\pi_{k+1}}(d_{\rho}^{\pi^{\star}})-V_{0}^{\pi_{k}}(d_{\rho}^{\pi^{\star}})+\lambda_k^{\top}\left(V_{g}^{\pi_{k+1}}(d_{\rho}^{\pi^{\star}})-V_{g}^{\pi_k}(d_{\rho}^{\pi^{\star}})\right)+\frac{3-\gamma}{1-\gamma}\|A_{L,k}^{\pi_k}-\widehat{A}_{L,k}\|.
    \label{lem:bounded_Lagrangian:eq2}
\end{align}

Combining \eqref{lem:bounded_Lagrangian:eq1} and \eqref{lem:bounded_Lagrangian:eq2}, we have
\begin{align*}
    V_{L,k}^{\pi^{\star}}(\rho)-V_{L,k}^{\pi_{k}}(\rho)
    &\leq \frac{1}{(1-\gamma)\alpha}\mathbb{E}_{s \sim d_{\rho}^{\pi^{\star}}}\left[D_{\text{KL}}(\pi^{\star}(\cdot\mid s)||\pi_k(\cdot\mid s))\right]-\frac{1}{(1-\gamma)\alpha}\mathbb{E}_{s \sim d_{\rho}^{\pi^{\star}}}\left[D_{\text{KL}}(\pi^{\star}(\cdot\mid s)||\pi_{k+1}(\cdot\mid s))\right]\notag\\
    &\hspace{20pt}+\frac{1}{1-\gamma}\left(V_{0}^{\pi_{k+1}}(d_{\rho}^{\pi^{\star}})-V_{0}^{\pi_{k}}(d_{\rho}^{\pi^{\star}})+\lambda_k^{\top}\left(V_{g}^{\pi_{k+1}}(s)-V_{g}^{\pi_k}(s)\right)+\frac{3-\gamma}{1-\gamma}\|A_{L,k}^{\pi_k}-\widehat{A}_{L,k}\|\right)\notag\\
    &\hspace{20pt}+\frac{1}{1-\gamma}\|A_{L,k}^{\pi_k}-\widehat{A}_{L,k}\|\notag\\
    &\leq \frac{1}{(1-\gamma)\alpha}\mathbb{E}_{s \sim d_{\rho}^{\pi^{\star}}}\left[D_{\text{KL}}(\pi^{\star}(\cdot\mid s)||\pi_k(\cdot\mid s))\right]-\frac{1}{(1-\gamma)\alpha}\mathbb{E}_{s \sim d_{\rho}^{\pi^{\star}}}\left[D_{\text{KL}}(\pi^{\star}(\cdot\mid s)||\pi_{k+1}(\cdot\mid s))\right]\notag\\
    &\hspace{20pt}+\frac{1}{1-\gamma}\left(V_{L,k}^{\pi_{k+1}}(d_{\rho}^{\pi^{\star}})-V_{L,k}^{\pi_{k}}(d_{\rho}^{\pi^{\star}})\right)+\frac{4-2\gamma}{(1-\gamma)^2}\|A_{L,k}^{\pi_k}-\widehat{A}_{L,k}\|
\end{align*}

Taking the sum from $k=0$ to $k=K-1$,
\begin{align}
    \frac{1}{K}\sum_{k=0}^{K-1}\left(V_{L,k}^{\pi^{\star}}(\rho)-V_{L,k}^{\pi_{k}}(\rho)\right) &\leq \frac{1}{(1-\gamma)K\alpha}\sum_{k=0}^{K-1}\mathbb{E}_{s\sim d_{\rho}^{\pi^{\star}}}\left[D_{\text{KL}}(\pi^{\star}(\cdot\mid s)||\pi_k(\cdot\mid s))-D_{\text{KL}}(\pi^{\star}(\cdot\mid s)||\pi_{k+1}(\cdot\mid s))\right]\notag\\
    &\hspace{20pt}+\frac{1}{(1-\gamma)K}\sum_{k=0}^{K-1}\left(V_{L,k}^{\pi_{k+1}}(d_{\rho}^{\pi^{\star}})-V_{L,k}^{\pi_{k}}(d_{\rho}^{\pi^{\star}})\right)+\frac{4-2\gamma}{(1-\gamma)^2 K}\sum_{k=0}^{K-1}\|A_{L,k}^{\pi_k}-\widehat{A}_{L,k}\|\notag\\
    &\leq \frac{1}{(1-\gamma)K\alpha}\mathbb{E}_{s\sim d_{\rho}^{\pi^{\star}}}\left[D_{\text{KL}}(\pi^{\star}(\cdot\mid s)||\pi_0(\cdot\mid s))\right]\notag\\
    &\hspace{20pt}+\frac{1}{(1-\gamma)K}\sum_{k=0}^{K-1}\left(V_{L,k}^{\pi_{k+1}}(d_{\rho}^{\pi^{\star}})-V_{L,k}^{\pi_{k}}(d_{\rho}^{\pi^{\star}})\right)+\frac{4-2\gamma}{(1-\gamma)^2 K}\sum_{k=0}^{K-1}\|A_{L,k}^{\pi_k}-\widehat{A}_{L,k}\|.
    \label{lem:bounded_Lagrangian:eq3}
\end{align}

The second term on the right hand side can be treated as follows
\begin{align}
    &\frac{1}{(1-\gamma)K}\sum_{k=0}^{K-1}\left(V_{L,k}^{\pi_{k+1}}(d_{\rho}^{\pi^{\star}})-V_{L,k}^{\pi_{k}}(d_{\rho}^{\pi^{\star}})\right)\notag\\
    &\leq \frac{1}{(1-\gamma)K}\sum_{k=0}^{K-1}\left(V_{0}^{\pi_{k+1}}(d_{\rho}^{\pi^{\star}})-V_{0}^{\pi_{k}}(d_{\rho}^{\pi^{\star}})\right)+\frac{1}{(1-\gamma)K}\sum_{k=0}^{K-1}\lambda_k^{\top}\left(V_{g}^{\pi_{k+1}}(d_{\rho}^{\pi^{\star}})-V_{g}^{\pi_{k}}(d_{\rho}^{\pi^{\star}})\right)\notag\\
    &= \frac{V_{0}^{\pi_{K}}(d_{\rho}^{\pi^{\star}})}{(1-\gamma)K}+\frac{1}{(1-\gamma)K}\sum_{k=0}^{K-1}\left(\lambda_{k+1}^{\top}V_{g}^{\pi_{k+1}}(d_{\rho}^{\pi^{\star}})-\lambda_k^{\top}V_{g}^{\pi_{k}}(d_{\rho}^{\pi^{\star}})\right)+\frac{1}{(1-\gamma)K}\sum_{k=0}^{K-1}\left(\lambda_{k}^{\top}-\lambda_{k+1}^{\top}\right)V_{g}^{\pi_{k+1}}(d_{\rho}^{\pi^{\star}})\notag\\
    &= \frac{V_{0}^{\pi_{K}}(d_{\rho}^{\pi^{\star}})}{(1-\gamma)K}+\frac{1}{(1-\gamma)K}\sum_{k=0}^{K-1}\left(\lambda_{k+1}^{\top}V_{g}^{\pi_{k+1}}(d_{\rho}^{\pi^{\star}})-\lambda_k^{\top}V_{g}^{\pi_{k}}(d_{\rho}^{\pi^{\star}})\right)+\frac{1}{(1-\gamma)K}\sum_{k=0}^{K-1}\left(\lambda_{k}^{\top}-\lambda_{k+1}^{\top}\right)V_{g}^{\pi_{k+1}}(d_{\rho}^{\pi^{\star}})\notag\\
    &= \frac{V_{0}^{\pi_{K}}(d_{\rho}^{\pi^{\star}})}{(1-\gamma)K}+\frac{1}{(1-\gamma) K}\sum_{i=1}^{N}\lambda_{i,K}V_i^{\pi_K}(d_{\rho}^{\pi^{\star}})+\frac{1}{(1-\gamma)K}\sum_{k=0}^{K-1}\sum_{i=1}^{N}(\lambda_{i,k}-\lambda_{i,k+1})V_{i}^{\pi_{k+1}}(d_{\rho}^{\pi^{\star}}).
    \label{lem:bounded_Lagrangian:eq4}
\end{align}

Due to the bounded reward functions, we know that $V_{i}^{\pi}(\zeta)\leq\frac{1}{1-\gamma}$ for all $\pi$, $\zeta$, and $i=1,2,\cdots,N$. Also the projection operator in dual variable update \eqref{Alg:AC:dual_update} guarantees that $|\lambda_{i,k}|$ is upper bounded by $\frac{2}{\xi(1-\gamma)}$. Plugging the bounds on the value functions and $|\lambda_{i,k}-\lambda_{i,k+1}|$ in \eqref{eq:lambda_ik_bound} into \eqref{lem:bounded_Lagrangian:eq4}, we have
\begin{align}
    |\lambda_{i,k}-\lambda_{i,k+1}|\leq \frac{2\eta}{1-\gamma},\,\forall i=1,2,\cdot,N, \,\,k=0,1,\cdots,K-1.
\end{align}

\begin{align*}
    &\frac{1}{(1-\gamma)K}\sum_{k=0}^{K-1}\left(V_{L,k}^{\pi_{k+1}}(d_{\rho}^{\pi^{\star}})-V_{L,k}^{\pi_{k}}(d_{\rho}^{\pi^{\star}})\right)\notag\\
    &\leq \frac{V_{0}^{\pi_{K}}(d_{\rho}^{\pi^{\star}})}{(1-\gamma)K}+\frac{1}{(1-\gamma) K}\sum_{i=1}^{N}\lambda_{i,K}V_i^{\pi_K}(d_{\rho}^{\pi^{\star}})+\frac{1}{(1-\gamma)K}\sum_{k=0}^{K-1}\sum_{i=1}^{N}(\lambda_{i,k}-\lambda_{i,k+1})V_{i}^{\pi_{k+1}}(d_{\rho}^{\pi^{\star}})\notag\\
    &\leq \frac{1}{(1-\gamma)^2 K}+\frac{2N}{\xi(1-\gamma)^3 K}+\frac{2N\eta}{(1-\gamma)^3} \leq\frac{3N}{\xi(1-\gamma)^3 K}+\frac{2N\eta}{(1-\gamma)^3}.
\end{align*}

Combining this with \eqref{lem:bounded_Lagrangian:eq3} implies
\begin{align*}
    &\frac{1}{K}\sum_{k=0}^{K-1}\left(V_{L,k}^{\pi^{\star}}(\rho)-V_{L,k}^{\pi_{k}}(\rho)\right)\notag\\
    &\leq \frac{1}{(1-\gamma)K\alpha}\mathbb{E}_{s\sim d_{\rho}^{\pi^{\star}}}\left[D_{\text{KL}}(\pi^{\star}(\cdot\mid s)||\pi_0(\cdot\mid s))\right]+\frac{3N}{\xi(1-\gamma)^3 K}+\frac{2N\eta}{(1-\gamma)^3}+\frac{4-2\gamma}{(1-\gamma)^2 K}\sum_{k=0}^{K-1}\|A_{L,k}^{\pi_k}-\widehat{A}_{L,k}\|\notag\\
    &\leq \frac{\log|\Acal|}{(1-\gamma)K\alpha}+\frac{3N}{\xi(1-\gamma)^3 K}+\frac{2N\eta}{(1-\gamma)^3}+\frac{4-2\gamma}{(1-\gamma)^2 K}\sum_{k=0}^{K-1}\|A_{L,k}^{\pi_k}-\widehat{A}_{L,k}\|,
\end{align*}
where the last inequality uses the fact that for $D_{\text{KL}}(p_1||p_2)\leq\log|\Acal|$ for $p_1,p_2\in\Delta_{\Acal}$ if $p_2$ is a uniform distribution.

\qed

\subsection{Proof of Lemma \ref{lem:step_improvement}}

From the performance difference lemma in \eqref{eq:def_lemma_performancediff},
\begin{align*}
V_{0}^{\pi_{k+1}}(\zeta)-V_{0}^{\pi_{k}}(\zeta)&= \frac{1}{1-\gamma}\mathbb{E}_{s \sim d_{\zeta}^{\pi_{k+1}}, a \sim \pi_{k+1}(\cdot \mid s)}\left[A_{0}^{\pi_k}(s, a)\right] \notag\\
&= \frac{1}{1-\gamma}\mathbb{E}_{s \sim d_{\zeta}^{\pi_{k+1}}, a \sim \pi_{k+1}(\cdot \mid s)}\left[Q_{0}^{\pi_k}(s, a)\right] -\frac{1}{1-\gamma}\mathbb{E}_{s \sim d_{\zeta}^{\pi_{k+1}}}\left[V_{0}^{\pi_k}(s)\right] \notag\\
&=\frac{1}{1-\gamma}\mathbb{E}_{s \sim d_{\zeta}^{\pi_{k+1}}, a \sim \pi_{k+1}(\cdot \mid s)}\left[\widehat{Q}_{L,k}(s, a)\right]-\frac{\lambda_k^{\top}}{1-\gamma}\mathbb{E}_{s \sim d_{\zeta}^{\pi_{k+1}}, a \sim \pi_{k+1}(\cdot \mid s)}\left[\widehat{Q}_{g,k}(s, a)\right]\notag\\
&\hspace{20pt}+\frac{1}{1-\gamma}\mathbb{E}_{s \sim d_{\zeta}^{\pi_{k+1}}, a \sim \pi_{k+1}(\cdot \mid s)}\left[Q_{0}^{\pi_k}(s,a)-\widehat{Q}_{0,k}(s, a)\right]-\frac{1}{1-\gamma}\mathbb{E}_{s \sim d_{\zeta}^{\pi_{k+1}}}\left[V_{0}^{\pi_k}(s)\right],
\end{align*}
where the last equality simply adds and subtracts the same terms and uses the relation $\widehat{Q}_{L,k}(s,a)=\widehat{Q}_{0,k}(s,a)+\lambda_k^{\top}\widehat{Q}_{g,k}(s,a)$. Note that the actor update rule \eqref{Alg:AC:actor_update_pi} implies
\begin{align*}
    \widehat{Q}_{L,k}(s, a)=\frac{1}{\alpha}\log\left(\frac{\pi_{k+1}(a\mid s)}{\pi_k(a\mid s)}\widehat{Z}_k(s)\right).
\end{align*}
Using this in the equality above, we have
\begin{align*}
V_{0}^{\pi_{k+1}}(\zeta)-V_{0}^{\pi_{k}}(\zeta) &=\frac{1}{1-\gamma}\mathbb{E}_{s \sim d_{\zeta}^{\pi_{k+1}}, a \sim \pi_{k+1}(\cdot \mid s)}\left[\widehat{Q}_{L,k}(s, a)\right]-\frac{\lambda_k^{\top}}{1-\gamma}\mathbb{E}_{s \sim d_{\zeta}^{\pi_{k+1}}, a \sim \pi_{k+1}(\cdot \mid s)}\left[\widehat{Q}_{g,k}(s, a)\right]\notag\\
&\hspace{20pt}+\frac{1}{1-\gamma}\mathbb{E}_{s \sim d_{\zeta}^{\pi_{k+1}}, a \sim \pi_{k+1}(\cdot \mid s)}\left[Q_{0}^{\pi_k}(s,a)-\widehat{Q}_{0,k}(s, a)\right]-\frac{1}{1-\gamma}\mathbb{E}_{s \sim d_{\zeta}^{\pi_{k+1}}}\left[V_{0}^{\pi_k}(s)\right] \notag\\
&=\frac{1}{(1-\gamma)\alpha}\mathbb{E}_{s \sim d_{\zeta}^{\pi_{k+1}}, a \sim \pi_{k+1}(\cdot \mid s)}\left[\log\left(\frac{\pi_{k+1}(a\mid s)}{\pi_k(a\mid s)}\widehat{Z}_k(s)\right)\right]\notag\\
&\hspace{20pt}-\frac{\lambda_k^{\top}}{1-\gamma}\mathbb{E}_{s \sim d_{\zeta}^{\pi_{k+1}}, a \sim \pi_{k+1}(\cdot \mid s)}\left[\widehat{A}_{g,k}(s, a)\right]-\frac{\lambda_k^{\top}}{1-\gamma}\mathbb{E}_{s \sim d_{\zeta}^{\pi_{k+1}}}\left[\widehat{V}_{g,k}(s)\right]\notag\\
&\hspace{20pt}+\frac{1}{1-\gamma}\mathbb{E}_{s \sim d_{\zeta}^{\pi_{k+1}}, a \sim \pi_{k+1}(\cdot \mid s)}\left[Q_{0}^{\pi_k}(s,a)-\widehat{Q}_{0,k}(s, a)\right]-\frac{1}{1-\gamma}\mathbb{E}_{s \sim d_{\zeta}^{\pi_{k+1}}}\left[V_{0}^{\pi_k}(s)\right] \notag\\
&=\frac{1}{(1-\gamma)\alpha}\mathbb{E}_{s \sim d_{\zeta}^{\pi_{k+1}}}\left[D_{\text{KL}}(\pi_{k+1}(\cdot\mid s)||\pi_{k}(\cdot\mid s))\right]+\frac{1}{(1-\gamma)\alpha}\mathbb{E}_{s \sim d_{\zeta}^{\pi_{k+1}}}\left[\log\widehat{Z}_k(s)\right]\notag\\
&\hspace{20pt}-\frac{\lambda_k^{\top}}{1-\gamma}\mathbb{E}_{s \sim d_{\zeta}^{\pi_{k+1}}, a \sim \pi_{k+1}(\cdot \mid s)}\left[\widehat{A}_{g,k}(s, a)\right]-\frac{\lambda_k^{\top}}{1-\gamma}\mathbb{E}_{s \sim d_{\zeta}^{\pi_{k+1}}}\left[\widehat{V}_{g,k}(s)\right]\notag\\
&\hspace{20pt}+\frac{1}{1-\gamma}\mathbb{E}_{s \sim d_{\zeta}^{\pi_{k+1}}, a \sim \pi_{k+1}(\cdot \mid s)}\left[Q_{0}^{\pi_k}(s,a)-\widehat{Q}_{0,k}(s, a)\right]-\frac{1}{1-\gamma}\mathbb{E}_{s \sim d_{\zeta}^{\pi_{k+1}}}\left[\widehat{V}_{0,k}(s)\right]\notag\\
&\hspace{20pt}-\frac{1}{1-\gamma}\mathbb{E}_{s \sim d_{\zeta}^{\pi_{k+1}}}\left[V_{0}^{\pi_k}(s)-\widehat{V}_{0,k}(s)\right].
\end{align*}

Since $\widehat{V}_{0,k}+\lambda_k^{\top}\widehat{V}_{g,k}=\widehat{V}_{L,k}$, $Q_{0}^{\pi_k}(s,a)-V_{0}^{\pi_k}(s)=A_{0}^{\pi_k}(s,a)$, and $\widehat{Q}_{0,k}(s,a)-\widehat{V}_{0,k}(s,a)=\widehat{A}_{0,k}(s,a)$, the equality further simplifies
\begin{align*}
V_{0}^{\pi_{k+1}}(\zeta)-V_{0}^{\pi_{k}}(\zeta) &= \frac{1}{(1-\gamma)\alpha}\mathbb{E}_{s \sim d_{\zeta}^{\pi_{k+1}}}\left[D_{\text{KL}}(\pi_{k+1}(\cdot\mid s)||\pi_{k}(\cdot\mid s))\right]+\frac{1}{(1-\gamma)\alpha}\mathbb{E}_{s \sim d_{\zeta}^{\pi_{k+1}}}\left[\log\widehat{Z}_k(s)\right]\notag\\
&\hspace{20pt}-\frac{\lambda_k^{\top}}{1-\gamma}\mathbb{E}_{s \sim d_{\zeta}^{\pi_{k+1}}, a \sim \pi_{k+1}(\cdot \mid s)}\left[\widehat{A}_{g,k}(s, a)\right]-\frac{\lambda_k^{\top}}{1-\gamma}\mathbb{E}_{s \sim d_{\zeta}^{\pi_{k+1}}}\left[\widehat{V}_{g,k}(s)\right]\notag\\
&\hspace{20pt}+\frac{1}{1-\gamma}\mathbb{E}_{s \sim d_{\zeta}^{\pi_{k+1}}, a \sim \pi_{k+1}(\cdot \mid s)}\left[Q_{0}^{\pi_k}(s,a)-\widehat{Q}_{0,k}(s, a)\right]-\frac{1}{1-\gamma}\mathbb{E}_{s \sim d_{\zeta}^{\pi_{k+1}}}\left[\widehat{V}_{0,k}(s)\right]\notag\\
&\hspace{20pt}-\frac{1}{1-\gamma}\mathbb{E}_{s \sim d_{\zeta}^{\pi_{k+1}}}\left[V_{0}^{\pi_k}(s)-\widehat{V}_{0,k}(s)\right] \notag\\
&= \frac{1}{(1-\gamma)\alpha}\mathbb{E}_{s \sim d_{\zeta}^{\pi_{k+1}}}\left[D_{\text{KL}}(\pi_{k+1}(\cdot\mid s)||\pi_{k}(\cdot\mid s))\right]+\frac{1}{(1-\gamma)\alpha}\mathbb{E}_{s \sim d_{\zeta}^{\pi_{k+1}}}\left[\log\widehat{Z}_k(s)\right]\notag\\
&\hspace{20pt}-\frac{\lambda_k^{\top}}{1-\gamma}\mathbb{E}_{s \sim d_{\zeta}^{\pi_{k+1}}, a \sim \pi_{k+1}(\cdot \mid s)}\left[\widehat{A}_{g,k}(s, a)\right]-\frac{1}{1-\gamma}\mathbb{E}_{s \sim d_{\zeta}^{\pi_{k+1}}}\left[\widehat{V}_{L,k}(s)\right]\notag\\
&\hspace{20pt}+\frac{1}{1-\gamma}\mathbb{E}_{s \sim d_{\zeta}^{\pi_{k+1}}, a \sim \pi_{k+1}(\cdot \mid s)}\left[A_{0}^{\pi_k}(s,a)-\widehat{A}_{0,k}(s, a)\right]\notag\\
&\geq \frac{1}{(1-\gamma)\alpha}\mathbb{E}_{s \sim d_{\zeta}^{\pi_{k+1}}}\left[\log\widehat{Z}_k(s)\right]-\frac{1}{(1-\gamma)}\mathbb{E}_{s \sim d_{\zeta}^{\pi_{k+1}}}\left[\widehat{V}_{L,k}(s)\right]\notag\\
&\hspace{20pt}-\frac{\lambda_k^{\top}}{(1\hspace{-2pt}-\hspace{-2pt}\gamma)}\mathbb{E}_{s \sim d_{\zeta}^{\pi_{k+1}}\hspace{-3pt},a\sim\pi_{k+1}(\cdot\mid s)}\hspace{-2pt}\left[\hspace{-1pt}A_{g}^{\pi_k}(s)\hspace{-1pt}\right]\hspace{-2pt}+\hspace{-2pt}\frac{1}{(1\hspace{-2pt}-\hspace{-2pt}\gamma)}\mathbb{E}_{s \sim d_{\zeta}^{\pi_{k+1}}\hspace{-3pt},a\sim\pi_{k+1}(\cdot\mid s)}\hspace{-2pt}\left[\hspace{-2pt}A_{L,k}^{\pi_k}(s,a)\hspace{-2pt}-\hspace{-2pt}\widehat{A}_{L,k}(s,a)\right],
\end{align*}
where the last inequality adds and subtracts the term $\frac{\lambda_k^{\top}}{(1-\gamma)}\mathbb{E}_{s \sim d_{\zeta}^{\pi_{k+1}}\hspace{-3pt},a\sim\pi_{k+1}(\cdot\mid s)}\left[A_{g}^{\pi_k}(s)\right]$ and uses the non-negativity of the Kullback–Leibler divergence.

Re-arranging this inequality and recognizing $\mathbb{E}_{s \sim d_{\zeta}^{\pi_{k+1}}\hspace{-3pt},a\sim\pi_{k+1}(\cdot\mid s)}\left[A_{g}^{\pi_k}(s)\right]=V_{g}^{\pi_{k+1}}(\zeta)-V_{g}^{\pi_k}(\zeta)$ due to the performance difference lemma in \eqref{eq:def_lemma_performancediff}, we obtain
\begin{align}
&V_{0}^{\pi_{k+1}}(\zeta)-V_{0}^{\pi_{k}}(\zeta)+\lambda_k^{\top}\left(V_{g}^{\pi_{k+1}}(\zeta)-V_{g}^{\pi_k}(\zeta)\right)\notag\\ &\geq\frac{1}{(1-\gamma)\alpha}\mathbb{E}_{s \sim d_{\zeta}^{\pi_{k+1}}}\left[\log\widehat{Z}_k(s)\right]-\frac{1}{(1-\gamma)}\mathbb{E}_{s \sim d_{\zeta}^{\pi_{k+1}}}\left[\widehat{V}_{L,k}(s)\right]\notag\\
&\hspace{20pt}+\frac{1}{(1-\gamma)}\mathbb{E}_{s \sim d_{\zeta}^{\pi_{k+1}},a\sim\pi_{k+1}(\cdot\mid s)}\left[A_{L,k}^{\pi_k}(s,a)-\widehat{A}_{L,k}(s,a)\right]\notag\\
&\geq\frac{1}{(1-\gamma)\alpha}\mathbb{E}_{s \sim d_{\zeta}^{\pi_{k+1}}}\left[\log\widehat{Z}_k(s)\right]-\frac{1}{(1-\gamma)}\mathbb{E}_{s \sim d_{\zeta}^{\pi_{k+1}}}\left[\widehat{V}_{L,k}(s)\right]-\frac{1}{1-\gamma}\|A_{L,k}^{\pi_k}-\widehat{A}_{L,k}\|.
\label{lem:step_improvement:eq1}
\end{align}

For the first term in the inequality \eqref{lem:step_improvement:eq1}, we have
\begin{align*}
    \log\widehat{Z}_k(s)&=\log \left(\sum_{a'\in\Acal}\pi_k(a' \mid s) \exp \left(\alpha \widehat{Q}_{L,k}(s, a')\right)\right)\notag\\
    &\geq \sum_{a'\in\Acal}\pi_k(a' \mid s)\log \left( \exp \left(\alpha \widehat{Q}_{L,k}(s, a')\right)\right)\notag\\
    &= \alpha \sum_{a'\in\Acal}\pi_k(a' \mid s_k)\widehat{Q}_{L,k}(s, a')\notag\\
    &= \alpha\widehat{V}_{L,k}(s) + \alpha \sum_{a'\in\Acal}\pi_k(a' \mid s_k)\widehat{A}_{L,k}(s, a')\notag\\
    &= \alpha\widehat{V}_{L,k}(s) + \alpha \sum_{a'\in\Acal}\pi_k(a' \mid s_k)\left(\widehat{A}_{L,k}(s, a')-A_{L,k}^{\pi_k}(s,a')\right),
\end{align*}
where the last inequality uses the fact that $\sum_{a'\in\Acal}\pi_k(a' \mid s_k)A_{L,k}^{\pi_k}(s,a')=0$. This bound on $\log\widehat{Z}_k(s)$ implies that 
\begin{align*}
    &\frac{1}{(1-\gamma)\alpha}\mathbb{E}_{s \sim d_{\zeta}^{\pi_{k+1}}}\left[\log\widehat{Z}_k(s)\right]-\frac{1}{(1-\gamma)}\mathbb{E}_{s \sim d_{\zeta}^{\pi_{k+1}}}\left[\widehat{V}_{L,k}(s)\right]\notag\\
    &=\frac{1}{(1-\gamma)\alpha}\sum_{s\in\Scal}d_{\zeta}^{\pi_{k+1}}(s)\left(\log\widehat{Z}_k(s)-\alpha \widehat{V}_{L,k}(s)-\alpha\sum_{a'\in\Acal}\pi_k(a' \mid s_k)\left(\widehat{A}_{L,k}(s, a')-A_{L,k}^{\pi_k}(s,a')\right)\right)\notag\\
    &\hspace{20pt}+\frac{1}{1-\gamma}\sum_{s\in\Scal}d_{\zeta}^{\pi_{k+1}}(s)\sum_{a'\in\Acal}\pi_k(a' \mid s_k)\left(\widehat{A}_{L,k}(s, a')-A_{L,k}^{\pi_k}(s,a')\right)\notag\\
    &\geq \frac{1}{\alpha}\sum_{s\in\Scal}\zeta(s)\left(\log\widehat{Z}_k(s)-\alpha \widehat{V}_{L,k}(s)-\alpha\sum_{a'\in\Acal}\pi_k(a' \mid s_k)\left(\widehat{A}_{L,k}(s, a')-A_{L,k}^{\pi_k}(s,a')\right)\right)\notag\\
    &\hspace{20pt}+\frac{1}{1-\gamma}\sum_{s\in\Scal}d_{\zeta}^{\pi_{k+1}}(s)\sum_{a'\in\Acal}\pi_k(a' \mid s_k)\left(\widehat{A}_{L,k}(s, a')-A_{L,k}^{\pi_k}(s,a')\right)\notag\\
    &\geq \frac{1}{\alpha}\mathbb{E}_{s \sim \zeta}\left[\log\widehat{Z}_k(s)\right]-\mathbb{E}_{s \sim \zeta}\left[\widehat{V}_{L,k}(s)\right]-\sum_{s\in\Scal}\zeta(s)\sum_{a'\in\Acal}\pi_k(a' \mid s_k)\left(\widehat{A}_{L,k}(s, a')-A_{L,k}^{\pi_k}(s,a')\right)\notag\\
    &\hspace{20pt}+\frac{1}{1-\gamma}\sum_{s\in\Scal}d_{\zeta}^{\pi_{k+1}}(s)\sum_{a'\in\Acal}\pi_k(a' \mid s_k)\left(\widehat{A}_{L,k}(s, a')-A_{L,k}^{\pi_k}(s,a')\right)\notag\\
    &\geq \frac{1}{\alpha}\mathbb{E}_{s \sim \zeta}\left[\log\widehat{Z}_k(s)\right]-\mathbb{E}_{s \sim \zeta}\left[\widehat{V}_{L,k}(s)\right]-\frac{2-\gamma}{1-\gamma}\|\widehat{A}_{L,k}-A_{L,k}^{\pi_k}\|,
\end{align*}
where the first inequality uses $d_{\zeta}^{\pi_{k+1}} \geq(1-\gamma) \zeta$ elementwise.

Combining this inequality with \eqref{lem:step_improvement:eq1} yields
\begin{align*}
    V_{0}^{\pi_{k+1}}(\zeta)-V_{0}^{\pi_{k}}(\zeta)+\lambda_k^{\top}\left(V_{g}^{\pi_{k+1}}(\zeta)-V_{g}^{\pi_k}(\zeta)\right)\geq\frac{1}{\alpha}\mathbb{E}_{s \sim \zeta}\left[\log\widehat{Z}_k(s)\right]-\mathbb{E}_{s \sim \zeta}\left[\widehat{V}_{L,k}(s)\right]-\frac{3-\gamma}{1-\gamma}\|A_{L,k}^{\pi_k}-\widehat{A}_{L,k}\|.
\end{align*}

\qed

\subsection{Proof of Lemma \ref{lem:constraint_violation}}

We define the perturbation function for $\phi\in\mathbb{R}^{N}$
\begin{align*}
    P(\phi)\triangleq \max_{\pi} V_0^{\pi}(\rho)\quad\text{subject to }&\quad V_i^{\pi}(\rho)\geq b_i+\phi_i\quad\forall i=1,\cdots,N.
\end{align*}

Then, we have for any $\pi$
\begin{align}
    V_{L}^{\pi,\lambda^{\star}}(\rho) \leq \max_{\pi} V_{L}^{\pi,\lambda^{\star}}(\rho)=V_{D}^{\lambda^{\star}}(\rho)=V_{0}^{\pi^{\star}}(\rho)=P(0),
    \label{lem:constraint_violation:eq1}
\end{align}
where the first equality uses the definition of the dual function, and the second equality comes from the strong duality.

Given $\phi \in \mathbb{R}^{N}$, we let $\pi_{\phi}$ denote that policy that attains the optimum for $P(\phi)$, i.e.
\begin{align*}
    \pi_{\phi}=\argmax V_0^{\pi}(\rho) \quad\text{subject to}\quad V_i^{\pi}(\rho)\geq b_i+\phi_i\quad\forall i=1,\cdots,N.
\end{align*}

Then, from \eqref{lem:constraint_violation:eq1} we have
\begin{align*}
    P(0)-\langle\phi,\lambda^{\star}\rangle&\geq V_{L}^{\pi_{\phi},\lambda^{\star}}(\rho)-\langle\phi,\lambda^{\star}\rangle=V_{0}^{\pi_{\phi}}(\rho)+\langle\lambda^{\star},V_{g}^{\pi_{\phi}}(\rho)-b\rangle-\langle\phi,\lambda^{\star}\rangle\notag\\
    &=V_{0}^{\pi_{\phi}}(\rho)+\langle\lambda^{\star},V_{g}^{\pi_{\phi}}(\rho)-b-\phi\rangle=P(\phi).
\end{align*}

Now, plugging $\phi=\bar{\phi}\triangleq-[b-V_g^{\bar{\pi}}(\rho)]_{+}$ into this result, we have
\begin{align}
    \langle\bar{\phi},\lambda^{\star}\rangle\leq P(0)- P(\bar{\phi}).
    \label{lem:constraint_violation:eq2}
\end{align}

Also, since $\bar{\phi}$ is non-positive,
\begin{align}
    V_0^{\bar{\pi}}(\rho)\leq V_0^{\pi^{\star}}(\rho)=P(0)\leq P(\bar{\phi}).
    \label{lem:constraint_violation:eq3}
\end{align}

Combining \eqref{lem:constraint_violation:eq2} and \eqref{lem:constraint_violation:eq3} implies
\begin{align*}
    \langle\bar{\phi},\lambda^{\star}\rangle\leq P(0)-P(\bar{\phi}) \leq V_0^{\pi^{\star}}(\rho)-V_0^{\bar{\pi}}(\rho).
\end{align*}
Then, since $\lambda^{\star}$ is entry-wise non-negative and $\bar{\phi}$ is entry-wise non-positive,
\begin{align*}
    (C-\|\lambda^{\star}\|_{\infty})\|\bar{\phi}\|_{1} &\leq \langle \lambda^{\star}-C\1, \bar{\phi}\rangle \leq V_0^{\pi^{\star}}(\rho)-V_0^{\bar{\pi}}(\rho)+C\|\bar{\phi}\|_1 =V_{0}^{\pi^{\star}}(\rho)-V_{0}^{\bar{\pi}}(\rho)+C\sum_{i=1}^{N}[b_i-V_i^{\bar{\pi}}(\rho)]_{+}\leq\delta.
\end{align*}

Therefore,
\begin{align*}
    \sum_{i=1}^{N}[b_i-V_i^{\bar{\pi}}(\rho)]_{+}=\|\bar{\phi}\|_{1}\leq \frac{\delta}{C-\|\lambda^{\star}\|_{\infty}}\leq \frac{2\delta}{C}.
\end{align*}

\qed

\subsection{Proof of Lemma \ref{lem:hatpi_onestepdrift}}
We have
\begin{align*}
    \|\widehat{\pi}_{k+1}-\widehat{\pi}_{k}\|\leq \|\pi_{k+1}-\pi_{k}\|\leq \|\theta_{k+1}-\theta_{k}\|
\end{align*}
where the second inequality is obvious from the definition of the behavior policy in \eqref{Alg:AC:behaviorpolicy_update}, and the third inequality is due to the softmax function being Lipschitz with constant 1.

From the update rule \eqref{Alg:AC:actor_update}, we have
\begin{align*}
    \|\theta_{k+1}-\theta_{k}\|^2 &= \alpha^2\|\widehat{Q}_{0,k}+\lambda_k^{\top} \widehat{Q}_{g,k}\|^2\leq 2\alpha^2\|\widehat{Q}_{0,k}\|^2+2\alpha^2\|\lambda_k^{\top} \widehat{Q}_{g,k}\|^2\notag\\
    &\leq 2\alpha^2\|\widehat{Q}_{0,k}\|^2+2\alpha^2 N\sum_{i=1}^{N}|\lambda_i|\| \widehat{Q}_{i,k}\|^2 \leq 2\alpha^2\frac{|\Scal||\Acal|}{(1-\gamma)^2}+\frac{4\alpha^2 N}{\xi(1-\gamma)}\sum_{i=1}^{N}\frac{|\Scal||\Acal|}{(1-\gamma)^2}\notag\\
    &\leq \frac{6\alpha^2 N^2 |\Scal||\Acal|}{\xi(1-\gamma)^3}.
\end{align*}

Combining the two inequalities above, we have
\begin{align*}
    \|\widehat{\pi}_{k+1}-\widehat{\pi}_{k}\|\leq \frac{3(|\Scal||\Acal|)^{1/2} N \alpha}{\xi^{1/2}(1-\gamma)^{3/2}}.
\end{align*}

\qed

\subsection{Proof of Lemma \ref{lem:z_onestepdrift}}
\begin{align*}
    \|z_{i,k+1}-z_{i,k}\|^2&=\|\widehat{Q}_{i,k+1}-Q_i^{\widehat{\pi}_{k+1}}-\widehat{Q}_{i,k}+Q_i^{\widehat{\pi}_{k}}\|^2\notag\\
    &\leq 2\|\widehat{Q}_{i,k+1}-\widehat{Q}_{i,k}\|^2+2\|Q_i^{\widehat{\pi}_{k+1}}-Q_i^{\widehat{\pi}_{k}}\|^2\notag\\
    &\leq 2\beta^2(\frac{2}{1-\gamma}+1)^2+2\|Q_i^{\widehat{\pi}_{k+1}}-Q_i^{\widehat{\pi}_{k}}\|^2 \notag\\
    &\leq 2\beta^2(\frac{2}{1-\gamma}+1)^2+\frac{18\gamma^2 (|\Scal||\Acal|)^3 N^2 \alpha^2}{\xi(1-\gamma)^7}\notag\\
    &\leq \frac{36 (|\Scal||\Acal|)^3 N^2 \beta^2}{\xi(1-\gamma)^7},
\end{align*}
where the second inequality comes from the update \eqref{Alg:AC:critic_update}, the third inequality applies Lemma \ref{lem:Q_Lipschitz_pi} and \ref{lem:hatpi_onestepdrift}, and the last inequality combines terms using $\alpha\leq\beta$.

Taking the square root yields
\begin{align*}
    \|z_{i,k+1}-z_{i,k}\|&\leq \frac{6 (|\Scal||\Acal|)^{3/2} N \beta}{\xi^{1/2}(1-\gamma)^{7/2}}.
\end{align*}

\qed

\subsection{Proof of Lemma \ref{lem:bound_Gamma}}

We introduce the following auxiliary Markov chain of states and actions generated under the policy $\theta_{k-\tau}$ starting from $s_{k-\tau}$ as follows
\begin{align}
&{s}_{k-\tau} \stackrel{\widehat{\pi}_{k-\tau-1}}{\longrightarrow}  \widetilde{a}_{k-\tau}\stackrel{\Pcal}{\longrightarrow}  \widetilde{s}_{k-\tau+1} \stackrel{\widehat{\pi}_{k-\tau-1}}{\longrightarrow} a_{k-\tau+1}\stackrel{\Pcal}{\longrightarrow}  \cdots\stackrel{\Pcal}{\longrightarrow}\widetilde{s}_{k} \stackrel{\widehat{\pi}_{k-\tau-1}}{\longrightarrow}  \widetilde{a}_{k}\stackrel{\Pcal}{\longrightarrow}  \widetilde{s}_{k+1} \stackrel{\widehat{\pi}_{k-\tau-1}}{\longrightarrow}  \widetilde{a}_{k+1}.
\label{eq:def_auxiliaryMarkovChain}
\end{align}
In comparison, we recall that the original Markovian samples generated by our algorithm are
\begin{align*}
&{s}_{k-\tau} \stackrel{\widehat{\pi}_{k-\tau-1}}{\longrightarrow}  a_{k-\tau}\stackrel{\Pcal}{\longrightarrow}  s_{k-\tau+1} \stackrel{\widehat{\pi}_{k-\tau}}{\longrightarrow} a_{k-\tau+1}\stackrel{\Pcal}{\longrightarrow}  \cdots\stackrel{\Pcal}{\longrightarrow} s_{k} \stackrel{\widehat{\pi}_{k-1}}{\longrightarrow}  a_{k}\stackrel{\Pcal}{\longrightarrow}  s_{k+1} \stackrel{\widehat{\pi}_{k}}{\longrightarrow}  a_{k+1}.
\end{align*}

We can decompose the term of interest
\begin{align}
    \mathbb{E}[\Gamma_i(\widehat{\pi}_k,z_{i,k},O_k)]&=\mathbb{E}[\Gamma_i(\widehat{\pi}_k,z_{i,k},O_k)-\Gamma_i(\widehat{\pi}_{k-\tau},z_{i,k},O_k)]\notag\\
    &\hspace{20pt}+\mathbb{E}[\Gamma_i(\widehat{\pi}_{k-\tau},z_{i,k},O_k)-\Gamma_i(\widehat{\pi}_{k-\tau},z_{i,k-\tau},O_k)]\notag\\
    &\hspace{20pt}+\mathbb{E}[\Gamma_i(\widehat{\pi}_{k-\tau},z_{i,k-\tau},O_k)-\Gamma_i(\widehat{\pi}_{k-\tau},z_{i,k-\tau},\widetilde{O}_k)]\notag\\
    &\hspace{20pt}+\mathbb{E}[\Gamma_i(\widehat{\pi}_{k-\tau},z_{i,k-\tau},\widetilde{O}_k)]
    \label{lem:bound_Gamma_decomposition}
\end{align}

We introduce the following lemmas to bound each term of \eqref{lem:bound_Gamma_decomposition}.

\begin{lem}[\citet{khodadadian2021finite} Lemma 11]\label{lem:Gamma_bound:lem1}
We have for all $\pi_1$, $\pi_2$, $z$, $O$
\begin{align*}
    \Gamma_i(\pi_1,z,O)-\Gamma_i(\pi_2,z,O)\leq C_2\|\pi_1-\pi_2\|,
\end{align*}
where $C_2=\frac{2\sqrt{2}|\Scal|^{3/2}|\Acal|^{3/2}}{(1-\gamma)^3}+\frac{8|\Scal|^2|\Acal|^3}{(1-\gamma)^2}\left(\left\lceil\log _{\rho} m^{-1}\right\rceil+\frac{1}{1-\rho}+2\right)$.

\end{lem}

\begin{lem}[\citet{khodadadian2021finite} Lemma 12]\label{lem:Gamma_bound:lem2}
We have for all $\pi$, $z_1$, $z_2$, $O$
\begin{align*}
    \Gamma_i(\pi,z_1,O)-\Gamma_i(\pi,z_2,O)\leq C_3\|z_1-z_2\|,
\end{align*}
where $C_3=1+\frac{9\sqrt{2|\Scal||\Acal|}}{1-\gamma}$.

\end{lem}

\begin{lem}[\citet{khodadadian2021finite} Lemma 13]\label{lem:Gamma_bound:lem3}
Defining $\Fcal_{t-\tau}=\{s_{t-\tau},\widehat{\pi}_{k-\tau-1},z_{i,k-\tau}\}$, we have
\begin{align*}
    \mathbb{E}[\Gamma_i(\widehat{\pi}_{k-\tau-1},z_{i,k-\tau},O_k)-\Gamma_i(\widehat{\pi}_{k-\tau-1},z_{i,k-\tau},\widetilde{O}_k)\mid F_{k-\tau}]\leq C_4\sum_{t=k-\tau}^{k}\mathbb{E}[\|\widehat{\pi}_{t}-\widehat{\pi}_{k-\tau-1}\|\mid F_{k-\tau}],
\end{align*}
where $C_4=\frac{4|\Scal|^{3/2}|\Acal|^{3/2}}{1-\gamma}(1+\frac{3|\Scal||\Acal|}{1-\gamma})$.

\end{lem}

\begin{lem}[\citet{khodadadian2021finite} Lemma 14]\label{lem:Gamma_bound:lem4}
Defining $\Fcal_{t-\tau}=\{s_{t-\tau},\widehat{\pi}_{k-\tau-1},z_{i,k-\tau}\}$, we have
\begin{align*}
    \mathbb{E}[\Gamma_i(\widehat{\pi}_{k-\tau-1},z_{i,k-\tau},\widetilde{O}_k)\mid F_{k-\tau}]\leq C_4\alpha.
\end{align*}

\end{lem}

Lemma \ref{lem:Gamma_bound:lem1} implies that the first term of \eqref{lem:bound_Gamma_decomposition} can be bounded as
\begin{align}
    \mathbb{E}[\Gamma_i(\widehat{\pi}_k,z_{i,k},O_k)-\Gamma_i(\widehat{\pi}_{k-\tau},z_{i,k},O_k)]&\leq C_2\mathbb{E}[\|\widehat{\pi}_k-\widehat{\pi}_{k-\tau}\|]\leq C_2\sum_{t=k-\tau}^{k-1}\mathbb{E}[\|\widehat{\pi}_{t+1}-\widehat{\pi}_{t}\|] \notag\\
    &\leq C_2 \sum_{t=k-\tau}^{k-1} \frac{3(|\Scal||\Acal|)^{1/2} N \alpha}{\xi^{1/2}(1-\gamma)^{3/2}}\leq \frac{3C_2(|\Scal||\Acal|)^{1/2} N \alpha\tau}{\xi^{1/2}(1-\gamma)^{3/2}},
    \label{lem:bound_Gamma:term1}
\end{align}
where the third inequality uses Lemma \ref{lem:hatpi_onestepdrift}.

Similarly, applying Lemma \ref{lem:Gamma_bound:lem2}, we consider the second term of \eqref{lem:bound_Gamma_decomposition}
\begin{align}
    \mathbb{E}[\Gamma_i(\widehat{\pi}_{k-\tau},z_{i,k},O_k)-\Gamma_i(\widehat{\pi}_{k-\tau},z_{i,k-\tau},O_k)]&\leq C_3 \mathbb{E}[\|z_{i,k}-z_{i,k-\tau}\|]\notag\\
    &\hspace{-200pt}\leq C_3\sum_{t=k-\tau}^{k-1}\mathbb{E}[\|z_{i,t+1}-z_{i,t}\|]\leq C_3\sum_{t=k-\tau}^{k-1}\frac{6 (|\Scal||\Acal|)^{3/2} N \beta}{\xi^{1/2}(1-\gamma)^{7/2}}\leq \frac{6C_3 (|\Scal||\Acal|)^{3/2} N \beta\tau}{\xi^{1/2}(1-\gamma)^{7/2}}.
    \label{lem:bound_Gamma:term2}
\end{align}

Using Lemma \ref{lem:Gamma_bound:lem3} and taking the expectation over $\Fcal_{k-\tau}$, we have for the following bound on the third term of \eqref{lem:bound_Gamma_decomposition},
\begin{align}
    \mathbb{E}[\Gamma_i(\widehat{\pi}_{k-\tau},z_{i,k-\tau},O_k)-\Gamma_i(\widehat{\pi}_{k-\tau},z_{i,k-\tau},\widetilde{O}_k)]&\leq C_4 \sum_{t=k-\tau}^{k}\mathbb{E}[\|\widehat{\pi}_{t}-\widehat{\pi}_{k-\tau-1}\|]\notag\\
    &\hspace{-200pt}\leq C_4\sum_{t=k-\tau}^{k}\tau\frac{3(|\Scal||\Acal|)^{1/2} N \alpha}{\xi^{1/2}(1-\gamma)^{3/2}}=\frac{3C_4(|\Scal||\Acal|)^{1/2} N \alpha\tau^2}{\xi^{1/2}(1-\gamma)^{3/2}},
    \label{lem:bound_Gamma:term3}
\end{align}
where the second inequality follows from Lemma \ref{lem:hatpi_onestepdrift}.

The fourth term of \eqref{lem:bound_Gamma_decomposition} follows similarly from taking the expectation over $\Fcal_{k-\tau}$ on the result of Lemma \ref{lem:Gamma_bound:lem4}
\begin{align}
    \mathbb{E}[\Gamma_i(\widehat{\pi}_{k-\tau},z_{i,k-\tau},\widetilde{O}_k)]\leq C_4\alpha
    \label{lem:bound_Gamma:term4}
\end{align}

Using \eqref{lem:bound_Gamma:term1}-\eqref{lem:bound_Gamma:term4} in \eqref{lem:bound_Gamma_decomposition}, we have
\begin{align*}
    \mathbb{E}[\Gamma_i(\widehat{\pi}_k,z_{i,k},O_k)] &\leq \frac{3C_2(|\Scal||\Acal|)^{1/2} N \alpha\tau}{\xi^{1/2}(1-\gamma)^{3/2}}+ \frac{6C_3 (|\Scal||\Acal|)^{3/2} N \beta\tau}{\xi^{1/2}(1-\gamma)^{7/2}}+\frac{3C_4(|\Scal||\Acal|)^{1/2} N \alpha\tau^2}{\xi^{1/2}(1-\gamma)^{3/2}}+C_4\alpha\notag\\
    &\leq \frac{(3C_2+6C_3+4C_4) (|\Scal||\Acal|)^{3/2} N \beta\tau^2}{\xi^{1/2}(1-\gamma)^{7/2}},
\end{align*}
where in the last inequality we combine terms using $|\Scal||\Acal|\geq 1$, $N\geq 1$, $\xi\leq 1$, and $\alpha\leq\beta$.

\qed






\end{document}